\newif\ifconf
\conffalse
\conftrue 

\documentclass[11pt]{article}
\usepackage{amssymb}
\usepackage{amsthm}
\usepackage{amsmath}
\usepackage{fullpage}
\usepackage{epsfig}
\usepackage{verbatim}
\usepackage{hyperref}

\newcommand\remove[1]{}

\newcommand{\rnote}[1]{}
\newcommand{\jnote}[1]{}

\newcommand{\Lip}{\mathrm{Lip}}

\newcommand{\vol}{\mathrm{vol}}

\newcommand{\e}{\varepsilon}
\newcommand{\F}{\mathbb F}

\newcommand{\E}{\mathbb{E}}

\DeclareMathOperator{\diam}{diam}

\newcommand{\R}{\mathbb{R}}
\newcommand{\PR}{\mathcal{P}}
\newcommand{\T}{\mathcal{T}}
\newcommand{\ED}{ {\rm ED}}

\newtheorem{theorem}{Theorem}[section]
\newtheorem{lemma}[theorem]{Lemma}

\newtheorem{corollary}[theorem]{Corollary}

\newtheorem{definition}[theorem]{Definition}

\newtheorem{remark}{Remark}[section]

\date{}
\begin{document}

\title{Nonembeddability theorems via Fourier analysis}

\author{ Subhash Khot\\Georgia Institute of Technology \\
    {\tt khot@cc.gatech.edu}
    \and Assaf Naor\\Microsoft Research \\
    {\tt anaor@microsoft.com}
    }

\maketitle

\begin{abstract}
Various new nonembeddability results (mainly into $L_1$) are
proved via Fourier analysis. In particular, it is shown that the
Edit Distance on $\{0,1\}^d$ has $L_1$ distortion $(\log
d)^{\frac12-o(1)}$. We also give new lower bounds on the $L_1$
distortion of flat tori, quotients of the discrete hypercube under
group actions, and the transportation cost (Earthmover) metric.
\end{abstract}

\section{Introduction}

The bi-Lipschitz theory of metric spaces has witnessed a surge of
activity in the past four decades. While the original motivation
for this type of investigation came from metric geometry and
Banach space theory, since the mid-1990s it has become
increasingly clear that understanding metric spaces in the
bi-Lipschitz category is intimately related to fundamental
algorithmic questions arising in theoretical computer science.
Despite the remarkable list of achievements of this field, which
includes the best known approximation algorithms for a a wide
range of NP hard problems, the bi-Lipschitz theory is still in its
infancy. In particular, there are very few known methods for
proving nonembeddability results. The purpose of this paper is
 to the develop a Fourier-analytic approach to
proving nonembeddability theorems. In doing so, we resolve several
problems, and shed new light on existing results. Additionally,
our work points toward several interesting directions for future
research, with emphasis on the study of the bi-Lipschitz structure
of quotients of metric spaces.

Let $(X,d_X)$ and $(Y,d_Y)$ be metric spaces. The Lipschitz
constant of a function $f:X\to Y$ is
$$
\|f\|_{\Lip}:=\sup_{\substack{x,y\in X\\x\ne y}}
\frac{d_Y(f(x),f(y))}{d_X(x,y)}.
$$
If $f$ is one-to-one then its distortion is defined as
$$
\mathrm{dist}(f):=\|f\|_{\Lip}\cdot \|f^{-1}\|_{\Lip}.
$$
If $f$ is not one-to-one then we set $\mathrm{dist}(f)=\infty$.
The least distortion with which $X$ can be embedded into $Y$ is
denoted $c_Y(X)$, namely
$$
c_Y(X):=\inf\{\mathrm{dist}(f):\ f:X\hookrightarrow Y\}.
$$

We are particularly interested in embeddings into $L_p$ spaces. In
this case we write $c_p(X)=c_{L_p}(X)$. The most studied type of
embeddings are  into Hilbert space, in which case the parameter
$c_2(X)$ is known as the Euclidean distortion of $X$. The
parameter $c_1(X)$, i.e. the least distortion required to embed
$X$ into $L_1$, is of great algorithmic significance, especially
in the study of cut problems in graphs. The Euclidean distortion
of a metric space $X$ is relatively well understood: it is enough
to understand the distortion of finite subsets of $X$, and for
finite metrics there is a simple semidefinite program which
computes their Euclidean distortion~\cite{LLR95}. Embeddings into
$L_1$ are much more mysterious (see~\cite{Linial02}), and there
are very few known methods to bound $c_1(X)$ from below.

The present paper contains several new nonembeddability results,
which we now describe. The common theme  is that our proofs are
based on analytic methods, most notably Fourier analysis on
$\{0,1\}^d$ and $\R^n$. We stress that this is not the first time
that nonembeddability results have drawn on techniques from
harmonic analysis. Indeed, the proofs of results
in~\cite{BMW86,Pisier86,NS02,KV05,MN05} all have a Fourier
analytic component.

\bigskip

\noindent {\bf Our results.}

\medskip

\noindent{\bf 1) Quotients of the discrete hypercube and
transportation cost.} A classical theorem of Banach states that
every separable Banach space is a quotient of $\ell_1$. More
precisely, for every separable Banach space $X$, there is a linear
subspace $Y\subseteq \ell_1$ such that $\ell_1/Y$ is linearly
isometric to $X$. This suggests that interesting ``bad examples"
of metric spaces can be obtained as {\em metric} quotients of the
Hamming cube. Roughly speaking, this says that we can obtain
interesting metrics (i.e. metrics that do not embed into nice
spaces, say $L_1$) by identifying points of the hypercube.
Quotients of metric spaces are a well studied concept
(see~\cite{Gro99,BH99} for an introduction, and~\cite{MN04} for a
discussion of quotients of finite metric spaces)- we refer the
reader to Section~\ref{section:quotients} for a precise definition
of this notion.

Motivated by this analogy, in Section~\ref{section:quotients} we
exhibit classes of quotients of the Hamming cube which do not
embed into $L_1$. A fundamental theorem of
Bourgain~\cite{Bourgain85} states that for every finite metric
space $X$, $c_1(X)\le c_2(X)=O(\log |X|)$. In~\cite{Bourgain85}
Bourgain used a counting argument to show that there exist
arbitrarily large metric spaces $X$ with
$c_2(X)=\Omega(\log|X|/\log\log |X|)$. In~\cite{LLR95,AR98} it was
shown that there exit arbitrarily large metric spaces $X$ with
$c_1(X)=\Omega(\log |X|)$ (namely $X$ can be taken to be a
constant degree expander). In Section~\ref{section:quotients} we
show that there exist simple $n$-point quotients of the Hamming
cube $\{0,1\}^d$ which incur distortion $\Omega(\log n)$ in any
$L_1$ embedding.
We also show that certain
quotients of the Hamming cube obtained from the action of a
transitive permutation group of the coordinates do not well-embed
into $L_1$. These results are proved via a flexible Fourier
analytic approach.

As an application of the results stated above we settle the
problem of the $L_1$ embeddability of the transportation cost
metric (also known as the Earthmover metric in the computer
vision/graphics literature) on the set of all probability measures
on $\{0,1\}^d$. Denoting by $\PR(\{0,1\}^d)$ the space of all
probability measures on the Hamming cube $\{0,1\}^d$, let
$\T_\rho(\sigma,\tau)$ denote the transportation cost distance
between $\sigma,\tau\in \PR(\{0,1\}^d)$, with respect to the cost
function induced by the Hamming metric $\rho$ (see
Section~\ref{section:transportation} for the definition). Such
metrics occur in various contexts in computer science: they are a
popular distance measure in graphics and vision~\cite{GRT00,IT03},
and they are used as LP relaxations for classification problems
such as $0$-extension and metric
labelling~\cite{CKNZ01,Cha02,AFHKTT04}. Transportation cost
metrics are also prevalent in several areas of analysis and PDEs
(see the book~\cite{Villani03} and the references therein).

Motivated by applications to nearest neighbor search (a.k.a.
similarity search in the vision literature), the problem of
embedding transportation cost metrics into $L_1$ attracted a lot
of attention in recent years (see~\cite{Cha02,IT03,jiriproblems}).
In~\cite{Cha02,IT03} it is shown that
$c_1(\PR(\{0,1\}^d),\T_\rho)=O(d)$. In
Section~\ref{section:transportation} we show that this bound is
optimal, i.e. $c_1(\PR(\{0,1\}^d),\T_\rho)=\Omega(d)$. From an
analytic perspective, Kantorovich duality (see~\cite{Villani03})
implies that $(\PR(\{0,1\}^d),\T_\rho)$ embeds isometrically into
$\Lip(\{0,1\}^d)^*$- the dual of the Banach space of all real
valued Lipschitz functions on the hypercube. A result of
Bourgain~\cite{Bourgain86-trees} implies that $\sup_{d\in \mathbb
N} c_1(\Lip(\{0,1\}^d)^*)=\infty$. Our result shows that in fact
$c_1(\Lip(\{0,1\}^d)^*)=\Omega(d)$, improving upon the lower bound
obtained in~\cite{Bourgain86-trees}.

\medskip

\noindent{\bf 2) Edit Distance does not embed into $L_1$.} Edit
Distance (also known as Levenstein distance~\cite{Lev65})  is a
metric defined on the set of all finite-length binary  strings,
which we denote $\{0,1\}^*$. This metric is best viewed as the
shortest path metric on the following infinite graph : Let $G$ be
a graph with set of vertices $\{0,1\}^*$, and $\{x,y\}$ is an edge
of the graph if the string $y$ can be obtained from string $x$ by
either deleting one character from $x$ or by inserting one
character into $x$. For strings $x, y$, denote the shortest path
distance in $G$ (i.e. the Edit Distance) between $x, y$ as
$\ED(x,y)$. In  words, $\ED(x,y)$ is the minimum number of edit
operations needed to transform $x$ into $y$.   Here we assume that
only insertion/deletion operations are allowed. Character
substitution can be simulated by a deletion followed by an
insertion. Similarly, one can {\it shift} a string by deleting its
first character and inserting it at the end.

Edit Distance is a very useful metric arising in several
applications, most notably in string and text comparison problems,
which are prevalent in computer science (e.g. compression and
pattern matching), computational biology, and web searching (see
the papers~\cite{MS00,CM02,ADGIR03,Ind04,BJKK04,OR05,CK05} and the
references therein, and the book~\cite{Gus97} for a discussion of
applications to computational biology).

Let $(\{0,1\}^d, \ED)$ denote the space $\{0,1\}^d$ with the Edit
Distance metric (inherited from the metric $\ED$ on $\{0,1\}^*$).
A well known problem, stated e.g. in~\cite{jiriproblems}, is
whether the space $(\{0,1\}^d,\ED)$ embeds into $L_1$ with
uniformly bounded distortion. Had this been true, it would have
had significant applications in computer science
(see~\cite{jiriproblems}). Most notably it would lead to
approximate nearest neighbor search algorithms under Edit
Distance, and to efficient algorithms for computing the Edit
Distance between two strings (both of these problems are being
solved, by computational biologists, every day, hundreds of times.
Getting a substantially faster algorithm for any of them would be
of great practical importance). In Section~\ref{section:edit} we
show that the $L_1$ embedding approach fails, by proving via
Fourier analytic methods that
\begin{eqnarray*}\label{eq:lower Edit}
c_1(\{0,1\}^d,\ED)\ge \frac{\sqrt{\log d}}{
2^{O\left(\sqrt{\log\log d\log\log\log d}\right)}}.
\end{eqnarray*}
The previous best known lower bound is due to~\cite{ADGIR03},
where it is shown that $c_1(\{0,1\}^d,\ED)\ge 3/2$. The best known
upper bound on $c_1(\{0,1\}^d,\ED)$ is due to~\cite{OR05}, where
it is proved that
$$
c_1(\{0,1\}^d,\ED)\le 2^{O\left(\sqrt{\log d\log\log d}\right)}.
$$

\remove{ The lower bound in~\eqref{eq:lower Edit} also holds true
for the {\em Block Edit Distance} on $\{0,1\}^d$. This distance is
defined analogously to the Edit Distance, except that we allow in
addition to the usual character insertion/deletion operations, the
operation ``substring move", which removes an arbitrary substring
and puts it anywhere else in the original string. Such metrics are
also studied in computational biology~\cite{MS00,CM02}, and
in~\cite{CM02} it is shown that the Block Edit Distance on
$\{0,1\}^d$ embeds into $L_1$ with distortion $O(\log d\log^* d)$.
No lower bound for the $L_1$ distortion of Block Edit Distance was
previously known.}

\medskip

\noindent{\bf 3) Flat tori can be highly non-Euclidean.} The Nash
embedding theorem~\cite{Nash54} states that any $n$-dimensional
Riemannian manifold is isometric to a Riemannian sub-manifold of
$\R^{2n}$. In the bi-Lipschitz category this is no longer the
case- it is easy to construct Riemannian manifolds (indeed, even
Riemannian surfaces) which do not embed bi-Lipschitzly even into
infinite dimensional Hilbert space. However, all the known
constructions were highly curved, and the possibility remained
that any manifold with zero curvature embeds bi-Lipschitzly into
$L_2$, with a uniform bound on the distortion. In
Section~\ref{section:torus} we show that this isn't the case:
there is an $n$-dimensional flat torus, i.e. $\R^n/\Lambda$ for
some lattice $\Lambda\subseteq \R^n$, equipped with the natural
Riemannian metric (whose sectional curvature is identically $0$),
such that $c_1(\R^n/\Lambda)=\Omega(\sqrt{n})$. This result
answers the question, posed by W. B. Johnson, whether a Lipschitz
quotient (in the sense of~\cite{BJLPS99}) of Hilbert space embeds
bi-Lipschitzly into Hilbert space. In~\cite{BJLPS99} it is shown
that a Banach space which is a Lipschitz quotient of a Hilbert
space is isomorphic to a Hilbert space. Johnson's question is
whether the condition that the quotient is a Banach space is
necessary. Since the natural quotient map $\pi: \R^n\to
\R^n/\Lambda$ is a Lipschitz quotient (see
Section~\ref{section:quotients}), the above example shows that
Lipschitz quotients of Hilbert space need not embed into Hilbert
space (indeed, they may not embed even into $L_1$). Our approach
is a variant of our study of quotient metrics in
Section~\ref{section:quotients}, and the proof is based on Fourier
analysis over $\R^n$, instead of discrete Fourier analysis over
$\{0,1\}^n$.

\bigskip

This paper is organized as follows. In
Section~\ref{section:discrete fourier} we present some background
and preliminary results on Fourier analysis on the Hamming cube.
In section~\ref{section:quotients} we investigate quotients of the
hypercube under group actions. In Section~\ref{section:edit} we
prove our lower bound on the $L_1$ distortion of Edit Distance,
and in Section~\ref{section:torus} we discuss the $L_1$ and $L_2$
embeddability of flat tori. We end with
Section~\ref{section:length}, which contains a brief discussion
which relates the notion of {\em length of metric spaces} (first
introduced by Schechtman~\cite{Sch82} in the context of the
concentration of measure phenomenon) to nonembeddability results.
This gives, in particular, new lower bounds on the Euclidean
distortion of various groups equipped with a group invariant
metric.

\section{Preliminaries on Fourier analysis on the
hypercube}\label{section:discrete fourier}

We start by introducing some notation concerning Fourier analysis
on the group $\F_2^d=\{0,1\}^d$. For $\e\in (0,1)$ we denote by
$\mu_\e$ the product $\e$-biased measure on $\F_2^d$, i.e. the
measure given by
$$
\forall\  x\in\F_2^d, \quad \mu_\e(\{x\})=\e^{\sum_{j=1}^d
x_j}(1-\e)^{d-\sum_{j=1}^d x_j}.
$$
For the sake of simplicity we write $\mu=\mu_{1/2}$. Given
$A\subseteq \{1,\ldots,d\}$ we define the Walsh function
$W_A:\F_2^d\to \R$ by
$$
W_A(x)=(-1)^{\sum_{j\in A}x_j}.
$$
Then $\{W_A:\ A\subseteq \{1,\ldots,d\}\}$ is an orthonormal basis
of $L_2(\F_2^d,\mu)$. In particular any $f:\F_2^d\to L_2$ has a
unique Fourier expansion
$$
f=\sum_{A\subseteq \{1,\ldots,d\}} \widehat f(A)W_A,
$$
where
$$
\widehat f(A)=\int_{\F_2^d} f(x)W_A(x)d\mu(x)\in L_2,
$$
and Parseval's identity reads as
$$
 \int_{\F_2^d} \|f(x)\|_2^2d\mu(x)=\sum_{A\subseteq \{1,\ldots, d\}}
\|\widehat f(A)\|_2^2.
$$
Let $e_j\in \F_2^d$ be the vector whose only non-zero coordinate
is the $j$th coordinate. We also write $e:=e_1+\ldots+e_d$ for the
all $1$s vector. The partial differentiation operator on
$L_2(\F_2^d)$ is defined by
$$
\partial_j f(x):=\frac{f(x+e_j)-f(x)}{2}.
$$
Since for every $A\subseteq \{1,\ldots, d\}$ we have that
$$
\partial_j W_A=\left\{\begin{array}{ll} -W_A & j\in A\\
0& j\notin A,\end{array}\right.
$$
we see that for every $f:\F_2^d\to \R$
\begin{eqnarray}\label{eq:as}
\sum_{j=1}^d \int_{\F_2^d}\partial _j
f(x)^2d\mu(x)=\sum_{A\subseteq \{1,\ldots,d\}}|A|\widehat f(A)^2.
\end{eqnarray}

In what follows we denote by $\rho$ the Hamming metric on
$\F_2^d$, namely for $x,y\in \F_2^d$,
$$
\rho(x,y):=|\{j\in \{1,\ldots, d\}:\ x_j\neq y_j\}|.
$$

Observe that for every $f:\F_2^d\to \R$,
$$
\int_{\F_2^d}
|f(x)-f(x+e)|^2d\mu(x)=\sum_{\substack{A\subseteq\{1,\ldots,d\}\\
|A|\equiv 1\mod 2}} 4\widehat f(A)^2\le 4\sum_{A\subseteq
\{1,\ldots,d\}}|A|\widehat f(A)^2= 4\sum_{j=1}^d
\int_{\F_2^d}[\partial _j f(x)]^2d\mu(x).
$$
This famous inequality, first proved by Enflo in~\cite{Enf69} via
a geometric argument, implies that $c_2(\F_2^d)\ge \sqrt{d}$.
Indeed, by integration we see that for every $f:\F_2^d\to L_2$,
$$
\int_{\F_2^d} \|f(x)-f(x+e)\|_2^2d\mu(x)\le  4\sum_{j=1}^d
\int_{\F_2^d}\|\partial _j f(x)\|_2^2d\mu(x).
$$
Thus, assuming that $f$ is invertible we see that
$$
\frac{d^2}{\|f^{-1}\|_{\Lip}^2}\le 4d\cdot
\left(\frac{\|f\|_{\Lip}}{2}\right)^2,
$$
i.e.
$$
\|f\|_{\Lip}\cdot \|f^{-1}\|_{\Lip}\ge \sqrt{d}.
$$
This Fourier-analytic approach to Enflo's theorem motivates the
ensuing arguments in this paper, since it turns out to be
remarkably flexible. For future reference we record here the basic
Poincar\'e inequality implied by the above reasoning:

\begin{lemma}\label{lem:hi fourier}
For every $f:\F_2^d \to L_2$,
$$
\int_{\F_2^d\times \F_2^d} \|f(x)-f(y)\|_2^2d\mu(x)d\mu(y)\le
\frac{2}{\min\{|A|:\ A\neq\emptyset,\ \widehat f(A)\neq
0\}}\sum_{j=1}^d \int_{\F_2^d}\|\partial _j f(x)\|_2^2d\mu(x).
$$
\end{lemma}

\begin{proof}
We simply observe that
$$
\int_{\F_2^d\times \F_2^d}
\|f(x)-f(y)\|_2^2d\mu(x)d\mu(y)=2\int_{\F_2^d} \|f(x)-\widehat
f(\emptyset)\|^2d\mu(x)=2\sum_{\emptyset \neq A\subseteq
\{1,\ldots, d\}} \|\widehat f(A)\|_2^2,
$$
and the required inequality follows from~\eqref{eq:as}.
\end{proof}

\section{Quotients of the hypercube}\label{section:quotients}

Let $(X,d_X)$ be a metric space. For $A,B\subseteq X$ the
Hausdorff distance between $A,B$ is defined as
\begin{eqnarray}\label{eq:hausdorff}
\mathcal H_X(A,B)=\sup\left\{\max\{d_X(a,B),d_X(b,A)\}:\ a\in A,\
b\in B\right\}.
\end{eqnarray}
Following~\cite{Gro99,BH99,MN04}, given a partition
$\mathcal{U}=\{U_1,\ldots,U_k\}$ of $X$, we define the {\em
quotient metric} induced by $X$ on $\mathcal U$, denoted
$X/\mathcal U$, as follows: assign to each $i,j\in\{1,\ldots,k\}$
the weight $w_{ij}=d_X(U_i,U_j)=\min_{x\in U_i, \ y\in
U_j}d_X(x,y)$, and let $d_{X/\mathcal U}(U_i,U_j)$ be the shortest
path distance between $i$ and $j$ in the weighted complete graph
on $\{1,\ldots,k\}$ in which the edge $\{i,j\}$ has weight
$w_{ij}$.

In the following lemma the right-hand inequality is an immediate
consequence of~\eqref{eq:hausdorff}, and the left-hand inequality
follows from the fact that the Hausdorff distance is a metric on
subsets of $X$.

\begin{lemma}\label{lem:no geo} Assume that $\mathcal
U=\{U_1,\ldots,U_k\}$ is a partition of a metric space $X$ such
that for every $i,j\in \{1,\ldots,k\}$, for every $x\in U_i$ there
exists $y\in U_j$ such that $d_X(x,y)=d_X(U_i,U_j)$. Then for
every $i,j\in\{1,\ldots,k\}$,
$$
d_{X/\mathcal U}(U_i,U_j)=\mathcal H_X(U_i,U_j)=d_X(U_i,U_j).
$$
\end{lemma}

A particular case of interest is when a group $G$ acts on $X$ by
isometries. In this case the orbit partition induced by $G$ on $X$
clearly satisfies the conditions of Lemma~\ref{lem:no geo},
implying that for all $x,y\in X$,
$$
d_{X/G}(Gx,Gy)=d_X(Gx,Gy),
$$
where we slightly abuse notation by letting $X/G$ be the quotient
of $X$ induced by the orbits of $G$. This is the only type of
quotients that we study in this paper. In particular,
Lemma~\ref{lem:no geo} implies  that the quotients we study here
are also {\em Lipschitz quotients} in the sense of~\cite{BJLPS99}
(see Section 6 in~\cite{MN04} for an explanation).

\remove{
\subsection{Warmup: the projective cube}

Consider the action of the cyclic group of order $2$, $C_2$, on
$\F_2^d$ given by the idempotent mapping $x\mapsto x+e$ (i.e. the
mapping obtained by flipping all the coordinates of $x\in
\F_2^d$). We call the quotient $\F_2^d/C_2$ the {\em projective
cube}. Consider a bi-Lipschitz mapping $f:\F_2^2/C_2\to L_2$ and
define $widetilde f: \F_2^d\to L_2$ via $\widetilde
f(x)=f(\{x,x+e\})$. Then $\widetilde f(x+x)=\widetilde f(x)$ for
every $x\in \F_2^d$, implying that for all $A\subseteq \{1,\ldots,
d\}$ with $|A|$ odd, $\widehat{(\widetilde f)}(A)=0$. By
Lemma~\ref{lem:hi fourier} we get that
$$
\frac{1}{\|f^{-1}\|_{\Lip}^2}\int_{\F_2^d\times \F_2^d}
\big[\rho_{\F_2^d/C_2}(\{x,x+e\},\{y,y+e\})\big]^2d\mu(x)d\mu(y)\le
\frac{d}{4}\cdot\|f\|_{\Lip}^2.
$$
}

\medskip

We will require the following lower bound on the average distance
in quotients of the hypercube.
\begin{lemma}\label{lem:average} Let $G$ be a group of isometries
acting on $\F_2^d$ with $2<|G|<2^d$. Then
$$
\int_{\F_2^d\times \F_2^d}
\rho_{\F_2^d/G}(Gx,Gy)d\mu(x)d\mu(y)=\Omega\left(
\frac{d-\log_2|G|}{1+\log_2\left(\frac{d}{d-\log_2|G|}\right)}\right).
$$
\end{lemma}

\begin{proof}
For every $t>0$,
\begin{eqnarray*}
\mu\times \mu\{x,y\in \F_2^d:\ \rho(Gx,Gy)\ge t\}\ge 1- \sum_{g\in
G}\mu\times\mu\{x,y\in \F_2^d:\ \rho(x,gy)\le t\}=
1-\frac{|G|}{2^{d}}\cdot \sum_{k\le t} \binom{d}{k}.
\end{eqnarray*}
We shall use the following (rough) bounds, which are a simple
consequence of Stirling's formula: For every $1/d< \delta\le 1/2$,
\begin{eqnarray}\label{eq:GV}
\frac{[\delta^\delta(1-\delta)^{1-\delta}]^{-d}}{6\sqrt{\delta
d}}\le \sum_{k\le \delta d} \binom{d}{k}\le 2\sqrt{\delta d}\cdot
[\delta^\delta(1-\delta)^{1-\delta}]^{-d}.
\end{eqnarray}
Thus, using Lemma~\ref{lem:no geo} we get that
\begin{eqnarray*}
\int_{\F_2^d\times \F_2^d} \rho_{\F_2^d/G}(Gx,Gy)d\mu(x)d\mu(y)\ge
\delta d\left(1-\frac{|G|}{2^d}2\sqrt{\delta d}\cdot
[\delta^\delta(1-\delta)^{1-\delta}]^{-d}\right).
\end{eqnarray*}
Choosing $\delta=\Theta\left(
\frac{d-\log_2|G|}{d+d\log_2\left(\frac{d}{d-\log_2|G|}\right)}\right)$
yields the required result.
\end{proof}

\subsection{A simple construction of $n$-point spaces with
$c_1=\Omega(\log n)$}

In what follows we refer to~\cite{MS77,BR98} for the necessary
background on coding theory. Let $C\subseteq \{0,1\}^d$ be a code,
i.e. a linear subspace of $\F_2^d$. Denote by $w(C)$ the minimum
Hamming weight of nonzero elements of $C$, i.e.
$$
w(C)=\min_{x\in C\setminus\{0\}} \|x\|_1.
$$
We also use the standard notation
$$
C^\perp:=\left\{x\in \F_2^d:\ \forall\ y\in C,\ \langle x,y\rangle
\equiv 0\mod 2\right\},$$
 where $\langle x,y\rangle :=\sum_{j=1}^n x_jy_j$.

\begin{lemma}\label{lem:hi frequency}
Assume that $f:\F^d_2\to L_2$ satisfies for every $x\in \F_2^d$
and $y\in C^\perp$, $f(x+y)=f(x)$. Then for every nonempty
$A\subseteq \{1,\ldots,d\}$ with $|A|<w(C)$, $\widehat f(A)=0$.
\end{lemma}
\begin{proof} Since $(C^\perp)^\perp=C$ (see~\cite{BR98}),
${\bf 1}_A\notin (C^\perp)^\perp$, implying that there exists
$v\in C^\perp$ such that $\langle {\bf 1}_A,v\rangle \equiv 1\mod
2$. Now,
\begin{eqnarray*}
\widehat f(A)&=&\int_{\F_2^n} f(x)W_A(x)d\mu(x)\\&=& \int_{\F_2^n}
f(x+v)W_A(x)d\mu(x)\\&=&\int_{\F_2^n}
f(x)W_A(x-v)d\mu(x)\\&=&(-1)^{\langle {\bf
1}_A,v\rangle}\int_{\F_2^n} f(x)W_A(x)d\mu(x)\\&=&-\widehat f(A).
\end{eqnarray*}
So $\widehat f(A)=0$.
\end{proof}

\begin{theorem}\label{thm:code} Let $C\subseteq \F_2^d$ be a code.
Then
$$
c_1(\F_2^d/C^\perp)=\Omega\left(w(C)\cdot\frac{\dim(C)}{d+d\log\left(\frac{d}{\dim(C)}\right)}\right).
$$
\end{theorem}

\begin{proof} Let $f:\F_2^d/C^\perp\to L_1$ be a bijection. Define
$\widetilde f:\F_2^d\to L_1$ by $\widetilde f(x)=f(x+C^\perp)$. It
is well known~\cite{DL97} that there exists a mapping $T:L_1\to
L_2$ such that for all $x,y\in L_1$,
$$\|T(x)-T(y)\|_2=\sqrt{\|x-y\|_1}.$$ Define $h:\F_2^d\to L_2$ by
$h=T\circ \widetilde f$. By Lemma~\ref{lem:hi frequency} and
Lemma~\ref{lem:hi fourier} we get that
\begin{eqnarray}\label{eq:upper code}
\int_{\F_2^d\times \F_2^d} \|\widetilde f(x)-\widetilde
f(y)\|_1d\mu(x)d\mu(y)\nonumber&=&\int_{\F_2^d\times \F_2^d}
\|h(x)-h(y)\|_2^2d\mu(x)d\mu(y)\\&\le&
\nonumber\frac{2}{w(C)}\sum_{j=1}^d \int_{\F_2^d}\|\partial _j
h(x)\|_2^2d\mu(x)\\\nonumber&=&\frac{2}{w(C)}\sum_{j=1}^d
\int_{\F_2^d}\|\partial _j \widetilde f(x)\|_1d\mu(x)\\
&\le& \frac{d}{w(C)}\|f\|_{\Lip}.
\end{eqnarray}
On the other hand, by Lemma~\ref{lem:average} we see that
\begin{eqnarray}\label{eq:lower code}
\int_{\F_2^d\times \F_2^d} \|\widetilde f(x)-\widetilde
f(y)\|_1d\mu(x)d\mu(y)\nonumber&=&\int_{\F_2^d\times \F_2^d} \|
f(x+C^\perp)-f(y+C^\perp)\|_1d\mu(x)d\mu(y)\\
&\ge& \frac{1}{\|f^{-1}\|_{\Lip}}\int_{\F_2^d\times
\F_2^d}\rho_{\F_2^d/C^\perp}(x+C^\perp,y+C^\perp)d\mu(x)d\mu(y)\nonumber\\
&=&\nonumber \Omega\left(
\frac{d-\log_2|C^\perp|}{1+\log_2\left(\frac{d}{d-\log_2|C^\perp|}\right)}\right)\cdot
\frac{1}{\|f^{-1}\|_{\Lip}}\\&=&
\Omega\left(\frac{\dim(C)}{1+\log\left(\frac{d}{\dim(C)}\right)}\right)\cdot
\frac{1}{\|f^{-1}\|_{\Lip}},
\end{eqnarray}
where we used the fact that $|C^\perp|=2^{d-\dim(C)}$.

Combining~\eqref{eq:upper code} and~\eqref{eq:lower code} yields
the required result.
\end{proof}

\begin{corollary}\label{coro:log n}
There exists arbitrarily large finite metric spaces $X$ for which
$c_1(X)=\Omega(\log |X|)$.
\end{corollary}

\begin{proof} Let $C\subseteq \{0,1\}^d$ be a code with
$ \dim(C)\ge \frac{d}{4}$ and $w(C)=\Omega(d)$. Such codes are
well known to exist (see~\cite{MS77}), and are easy to obtain via
the following greedy construction: fix $k\le d/4$ and let $V$ be a
$k$ dimensional subspace of $\F_2^d$ with $w(V)>\delta d$. Then
$V$ contains $2^k$ points. The number vectors $x\in \F_2^d$ with
$\|x+v\|_1\le \delta d$ for some $v\in V$ is at most
$2^k\sum_{\ell\le \delta d} \binom{d}{\ell}\le 2^{k+1}\sqrt{\delta
d}\cdot [\delta^\delta(1-\delta)^{1-\delta}]^{-d}$. It follows
that there exists $\delta=\Omega(1)$ such that for every $k\le
d/4$ there exists $x\in \F_2^d$ such that
$w(\mathrm{span}(V\cup\{x\}))>\delta d$, as required. Now, for $C$
as above, Theorem~\ref{thm:code} implies that
$$
c_1(\F_2^d/C^\perp)=\Omega(d)=\Omega(\log |\F_2^d/C^\perp|).
$$
\end{proof}

\begin{remark}\label{remark:extrapolation} {\em Using the Matou\v{s}ek's extrapolation lemma for
Poincar\'e inequalities~\cite{Mat97} (see also Lemma 5.5
in~\cite{BLMN05}), it is possible to prove that for a code $C$ as
in Corollary~\ref{coro:log n}, for every $p\ge 1$,
$c_p(\F_2^d/C^{\perp})\ge c(p) d$.}
\end{remark}

\subsection{The relation to transportation
cost}\label{section:transportation}

Given a finite metric space $(X,d)$ we denote by $\PR(X)$ the set
of all probability measures on $X$. For $\sigma,\tau\in \PR(X)$ we
define
$$
\Pi(\sigma,\tau) =\left\{\pi\in \PR(X\times X):\ \forall\, x\in X,
\ \int_{X}d\pi(x,y)=\sigma(x),\quad \mathrm{and}\quad
\int_{X}d\pi(y,x)=\tau(x) \right\}
$$
The optimal transportation cost (with respect to the metric $d$)
between $\sigma$ and $\tau$ is defined as
$$
\T_d(\sigma,\tau)=\inf_{\pi \in \Pi(\sigma,\tau)} \int_{X\times X}
d(x,y)d\pi(x,y).
$$
Given $A\subseteq X$ we denote by $\mu_A\in \PR(X)$ the uniform
probability measure on $A$. If $A,B\subseteq X$ have the same
cardinality then a straightforward extreme point argument
(see~\cite{Villani03}) shows that
$$
\T_d(\mu_A,\mu_B)=\inf\left\{ \int_A d(a,f(a))d\mu_A:\ f:A\to B\ \
\mathrm{is \ 1-1\ and\ onto}\right\}.
$$

\begin{lemma}\label{lem:group invariant} Let $G$ be a finite group,
equipped with a group invariant metric $d$ (i.e. $d(xg,yg)=d(x,y)$
for all $g,x,y\in G$). Then for every subgroup $H\subseteq G$ and
$x,y\in G$,
$$
d_{G/H}(xH,yH)=\T_d(\mu_{xH},\mu_{yH}).
$$
\end{lemma}
\begin{proof} For every bijection $f:xH\to yH$,
$$
\int_{xH} d(g,f(g))d\mu_{xH}(g)\ge d(xH,yH)= d_{G/H}(xH,yH).
$$

On the other hand, fix $h_1,h_2\in H$ such that
$d(xh_1,yh_2)=d(xH,yH)$. Then the mapping $f:xH\to yH$ given by
$f(g)=yh_2h_1^{-1}x^{-1}g$ satisfies for all $g\in xH$,
$d(g,f(g))=d(xh_1,yh_2)=d(xH,yH)$, implying the required result.
\end{proof}

\begin{corollary}\label{coro:transportation}
It follows from Corollary~\ref{coro:log n} and
Lemma~\ref{lem:group invariant} that $
c_1(\PR(\F_2^d),\T_\rho)=\Omega(d)$. This matches the upper bound
proved in~\cite{Cha02,IT03}. In fact, from
Remark~\ref{remark:extrapolation} we see that for all $p\ge 1$,
$c_p(\PR(\F_2^d),\T_\rho)\ge c(p)d.$
\end{corollary}

\begin{remark} {\em Let $\Lip(\F_2^d)$ be the Banach space of all functions $f:\F_2^d\to
\R$ satisfying $f(0)=0$, equipped with the Lipschitz norm $\|\cdot
\|_{\Lip}$. By Kantorovich duality (see~\cite{Villani03}),
$(\PR(\F_2^d),\T_\rho)$ is isometric to a subset of the dual space
$\Lip(\F_2^d)^*$. It follows that $c_1(\Lip(\F_2^d)^*)=\Omega(d)$.
As remarked in the introduction, the fact that $\sup_{d\in \mathbb
N} c_1(\Lip(\F_2^d)^*)=\infty$ was first proved by
Bourgain~\cite{Bourgain86-trees} using a different argument (which
yields a worse lower bound on the distortion).}
\end{remark}

\subsection{Actions of transitive permutation groups}

Let $G \le S_d$ be a subgroup of the symmetric group. Clearly $G$
acts by isometries on $\F_2^d$ via permutations of the
coordinates.

\begin{theorem}\label{thm:subgroup} Let $f:\F_2^d\to L_1$ be a
$G$-invariant function, where $G$ is transitive. Then
$$
\int_{\F_2^d\times\F_2^d} \|f(x)-f(y)\|_1d\mu(x)d\mu(y)\le
\frac{20}{\log d}\sum_{j=1}^d \int_{\F_2^d}
\|\partial_jf(x)\|_1d\mu(x).
$$
\end{theorem}

\begin{proof}
Let $A\subseteq\F_2^d$ be a $G$ invariant subset of the hypercube
and write $\mu(A)=p$. For $f={\bf 1}_A$ the required inequality
becomes:
\begin{eqnarray}\label{eq:boolean}
2p(1-p)\le \frac{10}{\log d}\sum_{j=1}^d I_j(A),
\end{eqnarray}
where $I_j(A)=\mu\{x\in \F_2^d:\ |\{x,x+e_j\}\cap A|=1\}|$ is the
{\em influence} of the $j$th variable on $A$. By~\cite{KKL88},
$\max_{1\le j\le d} I_j(A)\ge \frac{\log d}{5d}\cdot p(1-p)$. But,
since $A$ is invariant under the action of a transitive
permutation group, $I_j(A)$ is independent of $j$,
so~\eqref{eq:boolean} does indeed hold true.

In the general case let $f:\F_2^d\to L_1$ be a $G$ invariant
function. Denote by $\pi:\F_2^d\to \F_2^d/G$ the natural quotient
map, i.e. $\pi(x)=Gx$. Since $f$ is $G$-invariant, there is a
function $h: \F_2^d/G\to L_1$ such that $f=h\circ \pi$. By the
cut-cone representation of $L_1$ metrics (see~\cite{DL97}), there
are nonnegative weights $\{\lambda_A\}_{A\subseteq \F_2^d/G}$ such
that for every $x,y\in \F_2^d$, \begin{eqnarray*}
\|f(x)-f(y)\|_1&=&\|h(\pi(x))-h(\pi(y))\|_1\\&=&\sum_{A\subseteq
\F_2^d/G} \lambda_A|{\bf 1}_A(\pi(x))-{\bf
1}_A(\pi(y))|\\&=&\sum_{A\subseteq \F_2^d/G} \lambda_A|{\bf
1}_{\pi^{-1}(A)}(x)-{\bf 1}_{\pi^{-1}(A)}(y)|.
\end{eqnarray*}
Observe that for every $A\subseteq \F_2^d/G$,
$\pi^{-1}(A)\subseteq\F_2^d$ is $G$-invariant. Thus by the above
reasoning
\begin{eqnarray*}
\int_{\F_2^d\times\F_2^d} \|f(x)-f(y)\|_1d\mu(x)d\mu(y)&=&
\sum_{A\subseteq \F_2^d/G} \lambda_A\int_{\F_2^d\times\F_2^d}|{\bf
1}_{\pi^{-1}(A)}(x)-{\bf 1}_{\pi^{-1}(A)}(y)|d\mu(x)d\mu(y)\\
&\le& \sum_{A\subseteq \F_2^d/G} \lambda_A\cdot \frac{20}{\log
d}\sum_{j=1}^d\int_{\F_2^d} |\partial_j{\bf
1}_{\pi^{-1}(A)}(x)|d\mu(x)\\
&=& \frac{20}{\log d}\sum_{j=1}^d\int_{\F_2^d} \sum_{A\subseteq
\F_2^d/G} \lambda_A\left|\frac{{\bf 1}_{\pi^{-1}(A)}(x)-{\bf
1}_{\pi^{-1}(A)}(x+e_j)}{2}\right|d\mu(x)\\
&=&\frac{20}{\log d}\sum_{j=1}^d \int_{\F_2^d}
\|\partial_jf(x)\|_1d\mu(x).
\end{eqnarray*}
\end{proof}

We thus get many examples of spaces which do not well-embed into
$L_1$:

\begin{corollary}\label{coro:group} Let $G$ be a transitive
permutation group with $|G|<2^{\e d}$, for some $\e\in (0,1)$.
Then
$$
c_1(\F_2^d/G)\ge \Omega\left(\frac{(1-\e)}{1-\log(1-\e)}\cdot \log
d\right).
$$
\end{corollary}
\begin{proof}
This is a direct consequence of Theorem~\ref{thm:subgroup} and
lemma~\ref{lem:average}.
\end{proof}

\begin{remark}\label{rem:bourgain kalai}{\em
It is possible to obtain slightly stronger results analogous to
Corollary~\ref{coro:group} when we have additional information on
the structure of the group $G$. Indeed, in this case, in the proof
of Theorem~\ref{thm:subgroup}, one can use the results of Bourgain
and Kalai~\cite{BK97} on the influence of variables on group
invariant Boolean functions, instead of using~\cite{KKL88}.}
\end{remark}

\begin{remark} {\em We do not know if $(\F_2^d,\|\cdot\|_2)/G$ embeds
bi-Lipschitzly in Hilbert space with uniformly bounded distortion.
This seems to be unknown even in the case when $G$ is generated by
the cyclic shift of the coordinates. This problem is interesting
since if this space does embed into Hilbert space, then the
results of this section will yield an alternative approach to the
recent solution of the Goemans-Linial conjecture in~\cite{KV05}.}
\end{remark}

\section{Edit Distance does not embed into
$L_1$}\label{section:edit}

In this section we settle the $L_1$ embeddability problem of Edit
Distance negatively, by proving the following theorem:

\begin{theorem}\label{thm:edit} The following lower bound holds
true:
$$
c_1(\F_2^d,\ED)\ge \frac{\sqrt{\log d}}{ 2^{O\left(\sqrt{\log\log
d\log\log\log d}\right)}}.
$$
\end{theorem}

The following lemma is a useful way to prove $L_1$
nonembeddability results. The case $\delta=0$ of this lemma is due
to~\cite{LLR95}. Variants of the case $\delta>0$, which is the
case used in our proof of Theorem~\ref{thm:edit}, seem to be
folklore. We include here the formulation we need for the sake of
completeness (the main part of the proof below is a variant of the
proof of Lemma 3.6 in~\cite{NRS05}).

\begin{lemma}\label{lem:distributions} Fix $\alpha>0$ and
$0<\delta<\frac13$. Let $(X,d)$ be a finite metric space, $\sigma$
a probability measure on $X$, and $\tau$ a probability measure on
$X\times X$. Assume that for every $A\subseteq X$ with $\delta\le
\sigma(A)\le \frac23$ we have that $\tau(\{(x,y)\in X\times X:\
|\{x,y\}\cap A|=1\})\ge \alpha \sigma(A)$. Then,
$$
c_1(X)\ge \frac{\alpha}{2}\cdot\frac{\int_{X\times X} d(x,y)\,
d\sigma(x)d\sigma(y)-2\delta\diam(X)}{\int_{X\times X} d(x,y)\,
d\tau(x,y)}.
$$
\end{lemma}

\begin{proof} We claim that there exists a subset $Y\subseteq X$
with $\sigma(Y)\ge 1-\delta$ such that for every $f:Y\to L_1$,
\begin{eqnarray}\label{eq:subset}
\int_{Y\times Y} \|f(x)-f(y)\|_1d\sigma(x)d\sigma(y)\le
\frac{2}{\alpha}\int_{Y\times Y} \|f(x)-f(y)\|_1d\tau(x,y).
\end{eqnarray}
This will imply the required lower bound on $ c_1(X)$ since if
$f:X\to L_1$ is a bijection then
$$
\int_{Y\times Y} \|f(x)-f(y)\|_1d\tau(x,y)\le
\|f\|_{\Lip}\int_{X\times X} d(x,y)\, d\tau(x,y),
$$
while
\begin{eqnarray*}
\int_{Y\times Y} \|f(x)-f(y)\|_1d\sigma(x)d\sigma(y)&\ge&
\frac{1}{\|f^{-1}\|_{\Lip}}\left(\int_{X\times X} d(x,y)\,
d\sigma(x)d\sigma(y)-\right.\\&\phantom{\le}&\left.2\int_{ X\times
(X\setminus Y)} d(x,y)\,
d\sigma(x)d\sigma(y)\right)\\
&\ge& \frac{1}{\|f^{-1}\|_{\Lip}}\left(\int_{X\times X} d(x,y)\,
d\sigma(x)d\sigma(y)-2\delta\diam(X)\right).
\end{eqnarray*}

\bigskip

It remains to prove the existence of the required subset $Y$. For
simplicity we denote for every $A,B\subseteq X$,
$$
\{A,B\}=\Big\{(x,y)\in X\times X:\ \{x,y\}\cap
A\neq\emptyset\wedge \{x,y\}\cap B\neq \emptyset\Big\}.
$$
Define inductively disjoint subsets $\emptyset=
W_0,W_1,\ldots,W_k\subseteq X$ as follows. Having defined
$W_1,\ldots,W_i$, write $Y_i=\cup_{\ell=1}^i W_\ell$ and let
$W_{i+1}\subseteq X\setminus Y_i$ be an arbitrary nonempty subset
for which
$$\tau(\{W_{i+1},X\setminus (Y_i\cup W_{i+1})\})<\alpha \sigma(W_{i+1})\le
\frac{\alpha}{2} \sigma(X\setminus Y_i).$$ If no such $W_j$ exists
then this process terminates. We claim that $\sigma(Y_k)<\delta$.
Indeed, otherwise let $j$ be the first time at which
$\sigma(Y_{j})\ge\delta$. Observe that
$$
\sigma(Y_j)=\sigma(Y_{j-1})+\sigma(W_j)<\sigma(Y_{j-1})+\frac12\sigma(X\setminus
Y_{j-1})\le \frac{1+\delta}{2}\le \frac23.
$$
By our assumptions it follows that $\tau(\{Y_j,X\setminus
Y_j\})\ge \alpha \sigma(Y_j)$. But from the following simple
inclusion
$$
\{Y_j,X\setminus Y_j\}=\left\{\bigcup_{i=1}^j
W_i,X\setminus\bigcup_{i=1}^j W_i\right\}\subseteq \bigcup_{i=1}^j
\{W_i,X\setminus (Y_{i-1}\cup W_i)\}
$$
we deduce that
$$
0<\alpha\sigma(Y_j)\le \tau(\{Y_j,X\setminus
Y_j\})\le\sum_{i=1}^j\tau(\{W_i,X\setminus(Y_{i-1}\cup W_i)\})<
\sum_{i=1}^j \alpha \sigma(W_i)=\alpha\sigma(Y_j),
$$
a contradiction. Thus, taking $Y=Y_k$ we see that for every
$A\subseteq Y$ with $\sigma(A)\le\frac12$ we have
$\tau(\{A,Y\setminus A\})\ge \alpha \sigma(A)$. In other words,
\begin{eqnarray*}
\int_{Y\times Y} \|{\bf 1}_A(x)-{\bf
1}_A(y)\|_1d\tau(x,y)&=&\tau(\{A,Y\setminus A\})\\&\ge& \alpha
\sigma(A)\\&\ge& \alpha
\sigma(A)[\sigma(Y)-\sigma(A)]\\&=&\frac{\alpha}{2} \int_{Y\times
Y} \|f(x)-f(y)\|_1d\sigma(x)d\sigma(y),
\end{eqnarray*}
which implies~\eqref{eq:subset} by the cut cone representation of
$L_1$ metrics (as in the proof of Theorem~\ref{thm:subgroup}).
\end{proof}

In what follows we let $S$ denote the cyclic shift operator on
$\F_2^d$, namely
$$
S(x_1,\ldots,x_d)=(x_d,x_1,x_2,\ldots,x_{d-1}).
$$
We will also use the following remarkable theorem of
Bourgain~\cite{Bou02}. The exact dependence on the parameters
below is not stated explicitly in~\cite{Bou02}, but it follows
from the proof of~\cite{Bou02}. Since such quantitative estimates
are useful in several contexts, for future reference we reproduce
Bourgain's argument in Section~\ref{section:appendix}, while
tracking the bounds that he obtains.

\begin{theorem}[Bourgain's noise sensitivity
theorem~\cite{Bou02}]\label{thm:sensitive} Fix $\e,\delta\in
(0,1/10)$ and a Boolean function $f:\F_2^d\to \{-1,1\}$. Assume
that
$$
\sum_{A\subseteq \{1,\ldots,d\}}(1-\e)^{|A|}\widehat f(A)^2\ge
1-\delta.
$$
Then for every $\beta>0$ there exists a function $g:\F_2^d\to \R$
which depends on at most $\frac{1}{\e\beta}$ coordinates, such
that
$$
\int_{\F_2^d}\left(f(x)-g(x)\right)^2d\mu(x)\le
2^{c\sqrt{\log(1/\delta)\log\log(1/\e)}}\cdot
\left(\frac{\delta}{\sqrt{\e}}+4^{1/\e}\sqrt{\beta}\right).
$$
Here $c$ is a universal constant.
\end{theorem}

\begin{lemma}\label{lem:boolean edit} There exists a universal constant $C>0$ such that for every $\e\in (0,1/10)$, every integer $k\ge 10^{20/\e}$,
 and every $f:\F_2^d\to
\{-1,1\}$,
\begin{eqnarray}\label{eq:edit goal}
&&\!\!\!\!\!\!\!\!\!\!\!\!\!\!\!\!\int_{\F_2^d\times\F_2^d}
|f(x)-f(y)|d\mu(x)d\mu(y)-\nonumber
2^{-\sqrt{\log(1/\e)\log\log(1/\e)}}\\&\le&
\frac{2^{C\sqrt{\log(1/\e)\log\log(1/\e)}}}{k\sqrt{\e}}\cdot
\sum_{j=1}^{k}\int_{\F_2^d\times
\F_2^d}|f(x)-f(S^j(x)+y)|d\mu(x)d\mu_\e(y).
\end{eqnarray}
\end{lemma}
\begin{proof} Observe that for every $x,y\in \F_2^d$, $|f(x)-f(y)|=1-f(x)f(y)$. Thus for
every $j\in \{1,\ldots, d\}$
\begin{eqnarray}\label{eq:to beckner}
&&\!\!\!\!\!\!\!\!\!\!\!\!\!\int_{\F_2^d\times
\F_2^d}|f(x)-f(S^j(x)+y)|d\mu(x)d\mu_\e(y)=1-\int_{\F_2^d\times
\F_2^d}f(x)f(S^j(x)+y)d\mu(x)d\mu_\e(y)\nonumber\\
&=& 1-\int_{\F_2^d\times \F_2^d}\left(\sum_{A,B\subseteq
\{1,\ldots,d\}}\widehat f(A)\widehat f(B)
W_A(x)W_B(S^j(x)+y)\right)d\mu(x)d\mu_\e(y)\nonumber\\
&=& 1-\int_{\F_2^d\times \F_2^d}\left(\sum_{A,B\subseteq
\{1,\ldots,d\}}\widehat f(A)\widehat f(B)
W_A(x)W_{S^{-j}(B)}(x)W_B(y)\right)d\mu(x)d\mu_\e(y)\nonumber\\
&=& 1-\sum_{A\subseteq\{1,\ldots,d\}} (1-2\e)^{|A|} \widehat
f(A)\widehat f(S^{j}(A))\nonumber\\
&\ge& 1-\sum_{A\subseteq\{1,\ldots,d\}} (1-2\e)^{|A|} \widehat
f(A)^2,
\end{eqnarray}
where we used the Cauchy-Schwartz inequality and the facts that
for all $B\subseteq\{1,\ldots,d\}$ we have $\int_{\F_2^d}
W_B(y)d\mu_\e(y)=(1-2\e)^{|B|}$ and $\int_{\F_2^d}
W_AW_{S^{-j}(B)}d\mu=0$ when $B\neq S^{j}(A)$.
Averaging~\eqref{eq:to beckner} over $j=1,\ldots,k$ we see that
$$
\frac{1}{k}\sum_{j=1}^{k}\int_{\F_2^d\times
\F_2^d}|f(x)-f(S^j(x)+y)|d\mu(x)d\mu_\e(y)\ge
1-\sum_{A\subseteq\{1,\ldots,d\}} (1-2\e)^{|A|} \widehat f(A)^2.
$$
Thus, in order to prove~\eqref{eq:edit goal} we may assume that
\begin{eqnarray}\label{eq:for contradiction}
\sum_{A\subseteq\{1,\ldots,d\}} (1-2\e)^{|A|} \widehat f(A)^2\ge
1-\frac{\sqrt{\e}}{2^{4c\sqrt{\log(1/\delta)\log\log(1/\e)}}},
\end{eqnarray}
where $c>2$ is as in Theorem~\ref{thm:sensitive}. Now,  inequality
\eqref{eq:for contradiction}, together with
Theorem~\ref{thm:sensitive} (with $\beta=16^{-1/\e}$) implies that
there exists an integer $t\le 20^{1/\e}$, a function $g:\F_2^t\to
\{-1,1\}$, and indices $1\le i_1<i_2<\cdots <i_t\le d$ such that
if we extend $g$ to a function $\widetilde g:\F_2^d\to \{-1,1\}$
by setting
$$
\widetilde g(x_1,\ldots,x_d)=g(x_{i_1},x_{i_2},\ldots,x_{i_t}),
$$
then
$$
\int_{\F_2^d}|f(x)-\widetilde g(x)|d\mu(x)\le
\left(\int_{\F_2^d}\left(f(x)-\widetilde
g(x)\right)^2d\mu(x)\right)^2 \le
2^{-c\sqrt{\log(1/\e)\log\log(1/\e)}}.
$$
Write $I=\{i_1,\ldots,i_t\}$ and for $j\in \{1,\ldots,d\}$ define
$I+j=\{i_1+j \mod d,\ldots , i_t+j \mod d\}$. If $I\cap
(I+j)=\emptyset$ then we have the identity
\begin{eqnarray*}
\int_{\F_2^d\times \F_2^d}|\widetilde g(x)-\widetilde
g(S^j(x)+y)|d\mu(x)d\mu_\e(y)=\int_{ \F_2^d\times \F_2^d}
|\widetilde g(x)-\widetilde g(y)|d\mu(x)d\mu(y).
\end{eqnarray*}
Assume that $k\ge 2t^2$. In this case
$$|\{j\in \{1,\ldots, k\}:\
I\cap(I+j)=\emptyset\}|\ge k-t^2\ge \frac{k}{2}.$$ Thus
$$
\frac{1}{k}\sum_{j=1}^k \int_{\F_2^d\times \F_2^d}|\widetilde
g(x)-\widetilde g(S^j(x)+y)|d\mu(x)d\mu_\e(y)\ge \frac12
\int_{\F_2^d\times \F_2^d} |\widetilde g(x)-\widetilde
g(y)|d\mu(x)d\mu(y).
$$
It follows that
\begin{eqnarray*}
&&\!\!\!\!\!\!\!\!\!\!\!\!\!\!\!\!\frac{1}{k}\sum_{j=1}^k
\int_{\F_2^d\times \F_2^d}| f(x)- f(S^j(x)+y)|d\mu(x)d\mu_\e(y)\ge
\frac12 \int_{\F_2^d\times \F_2^d}
|f(x)-f(y)|d\mu(x)d\mu(y)-\\&\phantom{\le}&3\int_{\F_2^d}|f(x)-\widetilde
g(x)|d\mu(x)\\&\ge& \frac12 \int_{\F_2^d\times \F_2^d}
|f(x)-f(y)|d\mu(x)d\mu(y)-3\cdot
2^{-c\sqrt{\log(1/\e)\log\log(1/\e)}}.
\end{eqnarray*}
This completes the proof of~\eqref{eq:edit goal}.
\end{proof}

We require the following rough bound on the average edit distance
on $\F_2^d$. Such simple estimates have been previously obtained
by several authors, see for example Lemma 8 in~\cite{BEKMRRS03}.

\begin{lemma}\label{lem:edit average}
We have the following lower bound on the average Edit Distance on
$\F_2^d$:
$$
\int_{\F_2^d\times \F_2^d} \ED(x,y)d\mu(x)d\mu(y)\ge
\frac{d}{160}.
$$
\end{lemma}

\begin{proof} For every $x\in \F_2^d$ and every integer $r<d/2$,
$$
|\{y\in \F_2^d:\ \ED(x,y)=r\}|\le 2^r\binom{2d}{r}.
$$
This is best seen by observing that any sequence of $r$ insertions
or deletions can be executed in a sorted order, that is, the
indices of positions on which the operation is performed
increases. There are at most $\binom{2d}{r}$ ways to choose the
$r$ locations of these edit operations, and $2^r$ possible
insertion/deletion operations on these $r$ locations.

Now,
\begin{eqnarray*}
\mu\times \mu(\{(x,y)\in \F_2^d\times \F_2^d:\ \ED(x,y)>
d/16\})&\ge& 1-\frac{1}{2^d}\sum_{r\le d/16} 2^r\binom{2d}{r}\\
&\ge & 1-\frac{1}{2^d}\cdot 2^{d/8}\cdot2\sqrt{\frac{
d}{8}}\cdot\left[(1/16)^{1/16}(15/16)^{15/16}\right]^{-2d}\\&\ge&
\frac{1}{10}.
\end{eqnarray*}
\end{proof}

\begin{proof}[Proof of Theorem~\ref{thm:edit}] Let
$C$ be the constant in Lemma~\ref{lem:boolean edit}. Fix $\e
\in(0,1/10)$ such that $\e d>10^{20/\e}-1$, and an integer $\e
d\ge k\ge 10^{20/\e}$. Define a distribution $\tau$ on
$\F_2^d\times \F_2^d$ as follows: pick a pair $(x,y)\in
\F_2^d\times \F_2^d$ according to the measure $\mu\times \mu_\e$,
pick $j\in \{1,\ldots,k\}$ uniformly at random, and consider the
random pair $(x,S^j(x)+y)$. This induces a probability
distribution $\tau$ on $\F_2^d\times \F_2^d$. Observe that
$$
\ED(x,S^{j}(x)+y)\le 2\rho(0,y)+2j\le 2\rho(0,y)+2\e d.
$$
Thus
\begin{eqnarray}\label{eq:small dist}
\int_{\F_2^d\times \F_2^d}
\ED(x,y)d\tau(x,y)&=&\frac{1}{k}\sum_{j=1}^k \int_{\F_2^d\times
\F_2^d} \ED(x,S^{j}(x)+y)d\tau(x,y)\nonumber\\&\le& 2\e
d+2\sum_{r=0}^d \binom{d}{r}r\e^r(1-\e)^{d-r}=4\e d.
\end{eqnarray}
Lemma~\ref{lem:boolean edit} implies that for every $A\subseteq
\F_2^d$,
\begin{multline*}
\frac{2^{O\left(\sqrt{\log(1/\e)\log\log(1/\e)}\right)}}{\sqrt{\e}}\cdot
\tau(\{(x,y)\in\F_2^d\times \F_2^d:\ |\{x,y\}\cap A|=1\})\\ \ge
2\mu(A)[1-\mu(A)]- 3\cdot 2^{-\sqrt{\log(1/\e)\log\log(1/\e)}}.
\end{multline*}
Thus, the conditions of Lemma~\ref{lem:distributions} hold true
with the parameters $\delta=3\cdot
2^{-\sqrt{\log(1/\e)\log\log(1/\e)}}$ and
$\alpha=\frac{2^{O\left(\sqrt{\log(1/\e)\log\log(1/\e)}\right)}}{\sqrt{\e}}$.
Hence by~\eqref{eq:small dist} and Lemma~\ref{lem:edit average}
\begin{eqnarray*}
c_1(\F_2^d,\ED)\ge
\frac{\sqrt{\e}}{2^{O\left(\sqrt{\log(1/\e)\log\log(1/\e)}\right)}}\cdot\frac{\frac{d}{80}-6\cdot
2^{-\sqrt{\log(1/\e)\log\log(1/\e)}} \cdot 2d}{4\e d}.
\end{eqnarray*}
This implies the required result when we choose $\e\approx
\frac{1}{\log d}$.
\end{proof}

\remove{
\begin{remark} The same application of Lemma~\ref{lem:boolean
edit} yields the same lower bound for the $L_1$ distortion of
Block Edit Distance on $\F_2^d$. The only difference is that it is
necessary to replace Lemma~\ref{lem:edit average} by a lemma
stating that the average Block Edit Distance between two random
strings in $\F_2^d$ is $\Omega(d)$. This requires a
straightforward combinatorial argument which is slightly more
involved than the proof of Lemma~\ref{lem:edit average} (we omit
this argument here).
\end{remark}
}

\section{Flat tori which do not embed into
$L_1$}\label{section:torus}

Let $\Lambda\subseteq \R^n$ be a lattice in $\R^n$ of rank $n$.
The quotient space $\R^n/\Lambda$ is a Riemannian manifold
($n$-dimensional torus) whose curvature is identically zero.
Nevertheless, we show here that it is possible to construct
lattices $\Lambda$ such that $c_1(\R^n/\Lambda)=\Omega(\sqrt{n})$.
For a lattice $\Lambda\subseteq \R^n$ we denote its fundamental
parallelepiped by $P_\Lambda$. The dual lattice of $\Lambda$,
denoted $\Lambda^*$, is defined by
$$\Lambda^*=\{x\in \R^n:\ \forall\, y\in \Lambda,\
 \langle x,y\rangle\in \mathbb Z\}.$$
We shall use the following notation
$$
N(\Lambda)=\min_{x\in \Lambda\setminus \{0\}}\|x\|_2\quad
\mathrm{and}\quad r(\Lambda)=\max_{x\in \R^n}\min_{y\in \Lambda}
\|x-y\|_2.
$$
In words, $N(\Lambda)$ is the length of the shortest vector in
$\Lambda$, and $r(\Lambda)$ is the smallest $r$ such that balls of
radius $r$ centered at lattice points cover $\R^n$.

\begin{theorem}\label{thm:lattice}
Let $\Lambda\subseteq \R^n$ be a lattice. Then
$$
c_1(\R^n/\Lambda)=\Omega\left(\frac{N(\Lambda^*)}{r(\Lambda^*)}\cdot
\sqrt{n}\right).
$$
\end{theorem}

\begin{corollary}\label{coro:no torus}
Let $\Lambda$ be a lattice such that $\Lambda^*$ is almost
perfect, i.e. $N(\Lambda^*)=1$ and $r(\Lambda^*)\le 4$, say. Such
lattices are well known to exists (see~\cite{Zong02,Mar03}). Then
Theorem~\ref{thm:lattice} implies that
$c_1(\R^n/\Lambda)=\Omega(\sqrt{n})$. This is, in particular, an
example of a Riemannian manifold whose curvature is identically
zero which does not well-embed bi-Lipschitzly into $\ell_2$. This
fact should be contrasted with the Nash embedding
theorem~\cite{Nash54}, which says that any $n$-dimensional
Riemannian manifold is isometric to a Riemannian submanifold of
$\R^{2n}$.
\end{corollary}

\begin{remark} {\em Some restrictions on the Lattice $\Lambda$ should be
imposed in order to obtain a torus $\R^n/\Lambda$ which does not
embed into $\ell_2$. Indeed, the mapping $f:\R^n/\mathbb Z^n\to
\mathbb C^n$ defined by $f(x_1,\ldots,x_n)=(e^{2\pi i
x_1},\ldots,e^{2\pi i x_n})$ has distortion $O(1)$. We leave open
the interesting  problem of determining the value of
$c_1(\R^n/\Lambda)$ and $c_2(\R^n/\Lambda)$ as a function of
intrinsic geometric parameters of the lattice $\Lambda$. In
Theorem~\ref{thm:upper lattice} below we show that for every $n$
$$
L_n:=\sup\{c_2(\R^n/\Lambda):\ \Lambda\subseteq \R^n\ \mathrm{is\
a\ lattice}\}<\infty.
$$
Corollary~\ref{coro:no torus} shows that $L_n=\Omega(\sqrt{n})$,
while the upper bound obtained in Theorem~\ref{thm:upper lattice}
is $L_n=O(n^{3n/2})$. It would be of great interest to close the
large gap between these bounds.}
\end{remark}

\medskip

The proof of Theorem~\ref{thm:lattice} will be broken down into a
few lemmas. In what follows we fix a lattice $\Lambda\subseteq
\R^n$ and denote by $m$ the {\em normalized} Riemannian volume
measure on the torus $\R^n/\Lambda$. Given a function
$f:\R^n/\Lambda\to L_1$ we also think of $f$ as an
$\Lambda$-invariant function defined on $\R^n$. We refer
to~\cite{SW71} for the necessary background on Fourier analysis on
tori used in the ensuing arguments.

\begin{lemma}\label{lem:Cheeger lattice} Let $\gamma$ denote the
standard Gaussian measure on $\R^n$, i.e.
$d\gamma(x)=\frac{1}{(2\pi)^{n/2}} e^{-\|x\|_2^2/2}$. Then for
every continuous $f:\R^n/\Lambda\to L_1$,
\begin{eqnarray*}\label{eq:cheeger}
&&\!\!\!\!\!\!\!\!\!\!\!\!\!\!\int_{(\R^n/\Lambda)\times
(\R^n/\Lambda)} \|f(x)-f(y)\|_1dm(x)dm(y)\\&\le&
\frac{1}{1-e^{-2\pi^2[N(\Lambda^*)]^2}}\int_{\R^n}
\int_{\R^n/\Lambda} \|f(x)-f(x+y)\|_1dm(x)d\gamma(y).
\end{eqnarray*}
\end{lemma}

\begin{proof} By integration it is clearly enough to deal with the case of real-valued functions, i.e. $f:\R^n/\Lambda\to \R$.
Moreover, we claim that it suffices to prove the required
inequality when $f$ takes values in $\{0,1\}$. Indeed, assuming
the case of $f:\R^n/\Lambda\to \{0,1\}$, we pass to the general
case as follows:
\begin{eqnarray*}
&&\!\!\!\!\!\!\!\!\!\!\!\!\!\!\int_{(\R^/\Lambda)\times
(\R^n/\Lambda)} |f(x)-f(y)| dm(x)
dm(y)\\&=&\int_{(\R^n/\Lambda)\times (\R^n/\Lambda)}
\left(\int_{-\infty}^\infty|{\bf 1}_{(-\infty,t]}(f(x))-{\bf
1}_{(-\infty,t]}(f(y))|dt\right) dm(x)
dm(y)\\
&\le&\frac{1}{1-e^{-2\pi^2[N(\Lambda^*)]^2}}\int_{\R^n}
\int_{\R^n/\Lambda}\left(\int_{-\infty}^\infty|{\bf
1}_{(-\infty,t]}(f(x))-{\bf 1}_{(-\infty,t]}(f(x+y))|dt\right)
dm(x)d\gamma(y)\\
&=&\frac{1}{1-e^{-2\pi^2[N(\Lambda^*)]^2}}\int_{\R^n}
\int_{\R^n/\Lambda} |f(x)-f(x+y)|dm(x)d\gamma(y).
\end{eqnarray*}
So, it remains to prove the required inequality for a measurable
function $f:\R^n/\Lambda\to \{0,1\}$. The function $f$ can be
decomposed into a Fourier series indexed by the dual lattice
$\Lambda^*$:
$$
f(y)=\sum_{x\in \Lambda^*} \widehat f(x) e^{2\pi i\langle
x,y\rangle },
$$
where
$$
\widehat f (x)=\int_{\R^n/\Lambda} f(y)e^{-2\pi i\langle
x,y\rangle}dm(y).
$$
Using the fact that $|f(x)-f(x+y)|=f(x)+f(x+y)-2f(x)f(x+y)$ we get
from Parseval's identity that for every $y\in \R^n$
\begin{eqnarray*}
\int_{\R^n/\Lambda} |f(x)-f(x+y)|dm(x)&=& 2\widehat
f(0)-2\int_{\R^n/\Lambda}\left(\sum_{u,v\in \Lambda^*}\widehat
f(u)\widehat f(v)e^{2\pi i(\langle u,x\rangle+\langle
v,x+y\rangle)}\right)dm(x)\\
&=& 2\widehat f(0)-2\sum_{w\in \Lambda^*}e^{2\pi i \langle
w,y\rangle}|\widehat f(w)|^2.
\end{eqnarray*}
Integrating with respect to the Gaussian measure, and using the
identity $\int_{\R^n} e^{2\pi i \langle
w,y\rangle}d\gamma(y)=e^{-2\pi^2\|w\|_2^2}$, we get that
\begin{eqnarray}\label{eq:almost}
\int_{\R^n}\int_{\R^n/\Lambda} |f(x)-f(x+y)|dm(x)d\gamma(y)=
2\widehat f(0)[1-\widehat f(0)]-2\sum_{w\in
\Lambda^*\setminus\{0\}}e^{-2\pi^2\|w\|_2^2}|\widehat f(w)|^2.
\end{eqnarray}
On the other hand, since $f$ is Boolean function we have the
identities:
\begin{eqnarray}\label{eq:identities}
\int_{(\R^n/\Lambda)\times (\R^n/\Lambda)} |f(x)-f(y)| dm(x) dm(y)=
2\widehat f(0)[1-\widehat f(0)]=2\sum_{w\in \Lambda^*\setminus
\{0\}} |\widehat f(w)|^2.
\end{eqnarray}
Combining~\eqref{eq:almost} and~\eqref{eq:identities} we get
\begin{eqnarray*}
&&\!\!\!\!\!\!\!\!\!\!\!\!\!\!\!\!\int_{\R^n}\int_{\R^n/\Lambda}
|f(x)-f(x+y)|dm(x)d\gamma(y)= 2\sum_{w\in
\Lambda^*\setminus\{0\}}\left(1-e^{-2\pi^2\|w\|_2^2}\right)|\widehat
f(w)|^2\\&\ge& 2\left(1-e^{-2\pi^2[N(\Lambda^*)]^2}\right)\sum_{w\in
\Lambda^*\setminus\{0\}}|\widehat f(w)|^2\\
&=&
\left(1-e^{-2\pi^2[N(\Lambda^*)]^2}\right)\int_{(\R^n/\Lambda)\times
(\R^n/\Lambda)} |f(x)-f(y)| dm(x) dm(y).
\end{eqnarray*}
\end{proof}

\begin{lemma}\label{lem:average torus}
For every lattice $\Lambda\subseteq \R^n$,
$$
\int_{(\R^n/\Lambda)\times(\R^n/\Lambda) }
d_{\R^n/\Lambda}(x,y)dm(x)dm(y)\ge \frac{n}{16 r(\Lambda^*)}.
$$
\end{lemma}

\begin{proof} Let $V_\Lambda$ be the Voronoi cell of $\Lambda$
centered at $0$, i.e.
$$
V_\Lambda=\{x\in \R^n:\ \|x\|_2=d(x,\Lambda)\}.
$$
 Denote by $B_2^n$ the unit Euclidean ball
of $\R^n$ centered at $0$. Then by the definition of
$r(\Lambda^*)$ we have that $V_{\Lambda^*}\subseteq
r(\Lambda^*)B_2^n$. Hence $\vol(V_{\Lambda^*})\le
[r(\Lambda^*)]^n\vol(B_2^n)$. It is well known
(see~\cite{GL87,MG02,Mar03}) that
$$
\vol(V_\Lambda)\cdot
\vol(V_{\Lambda^*})=\vol(P_\Lambda)\cdot\vol(P_{\Lambda^*})=1.
$$
Thus
$$
\vol(V_\Lambda)\ge \frac{1}{[r(\Lambda^*)]^n\vol(B_2^n)}.
$$
It follows that
\begin{eqnarray*}
\int_{(\R^n/\Lambda)\times(\R^n/\Lambda) }
d_{\R^n/\Lambda}(x,y)dm(x)dm(y)&=&\frac{1}{\vol(V_\Lambda)}\int_{V_\Lambda}
\|x\|_2dx\\
&\ge& \frac{n}{8r(\Lambda^*)}\cdot \frac{\vol\left(\left\{x\in
V_\Lambda:\ \|x\|_2\ge
\frac{n}{8r(\Lambda^*)}\right\}\right)}{\vol(V_\Lambda)}\\
&\ge&
\frac{n}{8r(\Lambda^*)}\cdot\left(1-\left(\frac{n}{8r(\Lambda^*)}\right)^n\vol(B_2^n)\cdot
[r(\Lambda^*)]^n\vol(B_2^n)\right)\\
&\ge& \frac{n}{16r(\Lambda^*)}.
\end{eqnarray*}
\end{proof}


\begin{proof}[Proof of Theorem~\ref{thm:lattice}]
If $f:\R^n/\Lambda\to L_1$ is bi-Lipschitz then
\begin{eqnarray*}
\int_{\R^n} \int_{\R^n/\Lambda}
\|f(x)-f(x+y)\|_1dm(x)d\gamma(y)&\le& \|f\|_{\Lip}\cdot
\int_{\R^n} \int_{\R^n/\Lambda}
d_{\R^n/\Lambda}(x,x+y)dm(x)d\gamma(y)\\
&\le& \|f\|_{\Lip}\cdot \int_{\R^n} \|y\|_2d\gamma(y)\\
&\le& \|f\|_{\Lip}\cdot \sqrt{n}.
\end{eqnarray*}
On the other hand, using Lemma~\ref{lem:average torus} we see that
\begin{eqnarray*}
\int_{(\R^n/\Lambda)\times (\R^n/\Lambda)}
\|f(x)-f(y)\|_1dm(x)dm(y)&\ge&
\frac{1}{\|f^{-1}\|_{\Lip}}\int_{(\R^n/\Lambda)\times
(\R^n/\Lambda)}
d_{\R^n/\Lambda}(x,y)dm(x)dm(y)\\
 &\ge&
\frac{1}{\|f^{-1}\|_{\Lip}}\cdot\frac{n}{16r(\Lambda^*)},
\end{eqnarray*}
so by Lemma~\ref{lem:Cheeger lattice} we deduce that
$$
\|f\|_{\Lip}\cdot
\|f^{-1}\|_{\Lip}=\Omega\left(\frac{1-e^{-2\pi^2[N(\Lambda^*)]^2}}{r(\Lambda^*)}\cdot
\sqrt{n}\right).
$$
It follows that for every $t>0$,
$$
c_1(\R^n/\Lambda)=c_1(\R^n/(t\Lambda))=\Omega\left(\frac{1-e^{-2\pi^2[N((t\Lambda)^*)]^2}}{r((t\Lambda)^*)}\cdot
\sqrt{n}\right)=\Omega\left(\frac{1-e^{-2\pi^2[N(\Lambda^*)]^2/t^2}}{r(\Lambda^*)/t}\cdot
\sqrt{n}\right).
$$
Optimizing over $t$ yields the required result.
\end{proof}

\medskip

If one is interested only in bounding the Euclidean distortion of
$\R^n/\Lambda$, then the following lemma gives an alternative
proof of Theorem~\ref{thm:lattice} (in the case of embeddings into
$L_2$).

\begin{lemma}\label{lem:poincare lattice} For every continuous $f:\R^n/\Lambda\to
L_2$,
$$
\int_{(\R^n/\Lambda)\times (\R^n/\Lambda)}
\|f(x)-f(y)\|_2^2dm(x)dm(y)\le \frac{2}{\big[N(\Lambda^*)\big]^2}
\int_{\R^n/\Lambda} \|\nabla f(x)\|_2^2dm(x).
$$
\end{lemma}

\begin{proof}By Parseval's identity
\begin{eqnarray*}
 \int_{\R^n/\Lambda} \|\nabla
f(x)\|_2^2dm(x)&=&\sum_{j=1}^n
\int_{\R^n/\Lambda}\left(\frac{\partial f}{\partial x_j}(x)\right)^2dm(x)\\
&=& \sum_{x\in \Lambda^*} \|\widehat f (x)\|_2^2\cdot \|x\|_2^2 \\
&\ge& [N(\Lambda^*)]^2 \sum_{x\in \Lambda^*\setminus \{0\}} \|\widehat f (x)\|_2^2\\
&=& [N(\Lambda^*)]^2\int_{\R^n/\Lambda} \|f(x)-\widehat
f(0)\|_2^2dm(x)\\&=& \frac{[N(\Lambda^*)]^2}{2}
\int_{(\R^n/\Lambda)\times (\R^n/\Lambda)}
\|f(x)-f(y)\|_2^2dm(x)dm(y).
\end{eqnarray*}
\end{proof}

Lemma~\ref{lem:poincare lattice} yields the lower bound $
c_2(\R^n/\Lambda)=\Omega\left(\frac{N(\Lambda^*)}{r(\Lambda^*)}\cdot
\sqrt{n}\right)$ as follows. Let $f:\R^n/\Lambda\to L_2$ be a
bi-Lipschitz function. Since $L_2$ has the Radon-Nikodym property,
$f$ is differentiable almost everywhere (see~\cite{BL00}). Now, by
Lemma~\ref{lem:poincare lattice},
\begin{eqnarray*}
\int_{(\R^n/\Lambda)\times (\R^n/\Lambda)}
\|f(x)-f(y)\|_2^2dm(x)dm(y)&\le&
\frac{2}{\big[N(\Lambda^*)\big]^2} \sum_{j=1}^n\int_{\R^n/\Lambda}
\left\|\frac{\partial f}{\partial x_j}\right\|_2^2dm(x)\\&\le&
\frac{2}{\big[N(\Lambda^*)\big]^2}\cdot n\|f\|_{\Lip}^2.
\end{eqnarray*}
On the other hand, arguing as in the proof of
Theorem~\ref{thm:lattice}, we get
\begin{eqnarray*}
\int_{(\R^n/\Lambda)\times (\R^n/\Lambda)}
\|f(x)-f(y)\|_2^2dm(x)dm(y)&\ge&
\frac{1}{\|f^{-1}\|_{\Lip}^2}\int_{(\R^n/\Lambda)\times
(\R^n/\Lambda)}
d_{\R^n/\Lambda}(x,y)^2dm(x)dm(y)\\
&=&
\frac{1}{\|f^{-1}\|_{\Lip}^2}\cdot\Omega\left(\frac{n^2}{\big[r(\Lambda^*)\big]^2}\right).
\end{eqnarray*}
It follows that
$$
c_2(\R^n/\Lambda)=\Omega\left(\frac{N(\Lambda^*)}{r(\Lambda^*)}\cdot
\sqrt{n}\right).
$$

\bigskip

The following corollary of Lemma~\ref{lem:poincare lattice} will
not be used in the sequel, but we record it here for future
reference.

\begin{corollary}\label{lem:buser}
For every continuous $f:\R^n/\Lambda\to \R$,
$$
\int_{(\R^n/\Lambda)\times (\R^n/\Lambda)}
|f(x)-f(y)|dm(x)dm(y)\le \frac{2\sqrt{10}}{N(\Lambda^*)}
\int_{\R^n/\Lambda} \|\nabla f(x)\|_2dm(x).
$$
\end{corollary}

\begin{proof} Lemma~\ref{lem:poincare lattice} implies that
$\lambda_1(\R^n/\Lambda)\ge \big[N(\Lambda^*)\big]^2$, where
$\lambda_1(\R^n/\Lambda)$ is the smallest nonzero eigenvalue of
the Laplace-Beltrami operator on $\R^n/\Lambda$. Since
$\R^n/\Lambda$ has curvature $0$, an inequality of
Buser~\cite{Buser82} implies that $\lambda_1(\R^n/\Lambda)\le 10
\big[h(\R^n/\Lambda)\big]^2$, where $h(\R^n/\Lambda)$ is the
Cheeger constant of $\R^n/\Lambda$ (Buser's inequality can be
viewed as a reverse Cheeger inequality~\cite{Cheeger70} when the
Ricci curvature is bounded from below). Thus $h(\R^n/\Lambda)\ge
N(\Lambda^*)/\sqrt{10}$, which is precisely the required
inequality.
\end{proof}

\bigskip

We end this section by showing that there exists a constant
$D_n<\infty$ such that for any rank $n$ lattice $\Lambda\subseteq
\R^n$, $c_2(\R^n/\Lambda)\le D_n$.

\begin{lemma}\label{lem:korkin}  Every rank $n$ lattice
$\Lambda\subseteq \R^n$ has a basis (over $\mathbb Z^n$)
$x_1,\ldots,x_n$ such that for every $u_1,u_2,\ldots,u_n\in \R$,
\begin{eqnarray}\label{eq:2 estimate}
\frac{1}{n^{(3n-1)/2}}\cdot\left(\sum_{j=1}^n
u_j^2\|x_j\|_2^2\right)^{1/2}\le \left\|\sum_{j=1}^n u_j
x_j\right\|_2\le \sqrt{n}\cdot\left(\sum_{j=1}^n
u_j^2\|x_j\|_2^2\right)^{1/2}.
\end{eqnarray}
\end{lemma}
\remove{
\begin{proof} Let $\{x_1,\ldots,x_n\}$ be a basis of $\Lambda$.
Consider the {\em Gram-Schmidt orthogonalization} of this basis,
$\{y_1,\ldots,y_n\}$, which is defined recursively via $y_1=x_1$,
and
\begin{eqnarray}\label{eq:gram}
y_k=x_k-\sum_{j=1}^{k-1} \frac{\langle
x_k,y_j\rangle}{\|y_j\|_2^2}\cdot y_j.
\end{eqnarray}
Clearly $\|y_k\|_2\le \|x_k\|_2$ for all $k$. By a theorem
of~\cite{korkin} we can always chose the basis
$\{x_1,\ldots,x_n\}$ so that
$$
\prod_{k=1}^n\|x_k\|_2\le \sqrt{n}\cdot \prod_{k=1}^n\|y_j\|_2.
$$

Let $A$ be the matrix whose columns are the vectors
$\frac{x_1}{\|x_1\|_2},\ldots,\frac{x_n}{\|x_n\|_2}$ and let $B$
be the matrix whose columns are the vectors
$\frac{y_1}{\|y_1\|_2},\ldots,\frac{y_n}{\|y_n\|_2}$. Then $B$ is
an orthogonal matrix. By~\eqref{eq:gram} the matrix $BA$ is upper
triangular, with diagonal
$(\|y_1\|_2/\|x_1\|_2,\ldots,\|y_n\|_2/\|x_n\|_2)$. For an
$n\times n$ matrix $C$ we denote by $s_1(C)\ge s_2(C)\ge \cdots
s_n(C)> 0$ its singular values. Given any vector
$u=(u_1,\ldots,u_n)\in \R^n$ we have by the Cauchy-Schwartz
inequality that
$$
\|A u\|_2=\left\|\sum_{j=1}^n \frac{u_j}{\|x_j\|_2}\cdot
x_j\right\|_2\le \sum_{j=1}^n |u_j|\le \sqrt{n}\cdot \|u\|_2.
$$
This proves the righthand-side of~\eqref{eq:2 estimate}, and also
shows that $s_1(A)\le \sqrt{n}$. Since $B$ is orthogonal,
\end{proof}
}

\begin{proof}
Let $\{x_1,\ldots,x_n\}$ be a basis of $\Lambda$, and denote by
$A$ the matrix whose columns are the vectors
$\frac{x_1}{\|x_1\|_2},\ldots,\frac{x_n}{\|x_n\|_2}$. If we let
$\{x_1,\ldots,x_n\}$ be the Korkin-Zolotarev basis of $\Lambda$,
we can ensure that (see~\cite{LLS90}):
$$
|\det(A)|\ge \frac{1}{n^n}.
$$

Denote by $s_1(A)\ge s_2(A)\ge \cdots s_n(A)> 0$ the singular
values of $A$. Given a vector $u=(u_1,\ldots,u_n)\in \R^n$ we have
by the Cauchy-Schwartz inequality that
$$
\|A u\|_2=\left\|\sum_{j=1}^n \frac{u_j}{\|x_j\|_2}\cdot
x_j\right\|_2\le \sum_{j=1}^n |u_j|\le \sqrt{n}\cdot \|u\|_2.
$$
This proves the right-hand side of~\eqref{eq:2 estimate}, and also
shows that $s_1(A)\le \sqrt{n}$. Now
$$
\frac{1}{n^n}\le |\det(A)|=\prod_{j=1}^n s_j(A)\le s_1(A)\cdot
[s_n(A)]^{n-1}\le s_1(A)\cdot n^{(n-1)/2},
$$
i.e. $s_1(A)\ge n^{-(3n-1)/2}$. It follows that for every $u\in
R^n$, $\|Au\|_2\ge s_1(A)\|u\|_2\ge n^{-(3n-1)/2}\|u\|_2$, which
is precisely the left-hand side of~\eqref{eq:2 estimate}.
\end{proof}

\begin{theorem}\label{thm:upper lattice} Let $\Lambda\subseteq
\R^n$ be a lattice of rank $n$. Then $\R^n/\Lambda$ embeds into
$\R^{2n}$ with distortion $O(n^{3n/2})$.
\end{theorem}

\begin{proof} Let $\{x_1,\ldots,x_n\}$ be a basis as in
Lemma~\ref{lem:korkin}. Define $f:\R^n\to \mathbb C^n$ by
$$
f\left(\sum_{j=1}^n a_j x_j\right)=\left(\|x_1\|_2e^{2\pi i
a_1},\ldots,\|x_n\|_2e^{2\pi i a_n}\right).
$$
Since $f$ is $\Lambda$-invariant, we may think of it a a function
defined on the torus $\R^n/\Lambda$. For every $t\in \R$ let
$m(t)$ be the unique integer such that $t-m(t)\in [-1/2,1/2)$.
Given $u,v\in \R^n$,
\begin{eqnarray*}
\left\|f\left(\sum_{j=1}^n u_j x_j\right)-f\left(\sum_{j=1}^n v_j
x_j\right)\right\|_2^2&=&\sum_{j=1}^n \left| e^{2\pi
i(u_j-v_j)}-1\right|^2\cdot
\|x_j\|_2^2\\&=&2\sum_{j=1}^n\left[1-\cos(2\pi
(u_j-v_j))\right]\cdot \|x_j\|_2^2.
\end{eqnarray*}
Since for every $t\in \R$,
$$
\frac{[t-m(t)]^2}{12}\le 1-\cos(2\pi t)\le \frac{[t-m(t)]^2}{2},
$$
we get that
$$
\left\|f\left(\sum_{j=1}^n u_j x_j\right)-f\left(\sum_{j=1}^n v_j
x_j\right)\right\|_2^2=\Theta\left(\sum_{j=1}^n
[u_j-v_j-m(u_j-v_j)]^2\|x_j\|_2^2\right).
$$
On the other hand, by~\eqref{eq:2 estimate},
\begin{eqnarray*}
d_{\R^n/\Lambda}\left(\sum_{j=1}^n u_j x_j,\sum_{j=1}^n v_j
x_j\right)&=&d_{\R^n}(\sum_{j=1}^n (u_j-v_j) x_j,\Lambda)\\&\le&
\left\|\sum_{j=1}^n (u_j-v_j)
x_j-\sum_{j=1}^nm(u_j-v_j)x_j\right\|_2\\&\le&
\sqrt{n}\left(\sum_{j=1}^n
[u_j-v_j-m(u_j-v_j)]^2\|x_j\|_2^2\right)^{1/2}.
\end{eqnarray*}
In the reverse direction, let $m_1,\ldots,m_n\in \mathbb Z$ be
such that $\sum_{j=1}^n m_j x_j\in \Lambda$ is a closest lattice
point to $u-v$. Then
\begin{eqnarray*}
d_{\R^n/\Lambda}\left(\sum_{j=1}^n u_j x_j,\sum_{j=1}^n v_j
x_j\right)&=& \left\|\sum_{j=1}^n [
u_j-v_j-m_j]x_j\right\|_2\\
&\ge& \frac{1}{n^{(3n-1)/2}}\left(\sum_{j=1}^n [
u_j-v_j-m_j]^2\cdot\|x_j|\|_2^2\right)^{1/2}\\
&\ge& \frac{1}{n^{(3n-1)/2}}\left(\sum_{j=1}^n [
u_j-v_j-m(u_j-v_j)]^2\cdot\|x_j\|_2^2\right)^{1/2}.
\end{eqnarray*}
It follows that $f$ has distortion $O(n^{3n/2})$.
\end{proof}

\section{Length of metric spaces}\label{section:length}

The following definition, due to G. Schechtman~\cite{Sch82}, plays
an important role in the study of the concentration of measure
phenomenon and Levy families~\cite{MS86,Led01}.

\begin{definition}\label{def:length}
Let $(X,d)$ be a finite metric space. The length of $(X,d)$,
denoted $\ell(X,d)$ is the least constant $\ell$ such that there
exists a sequence of partitions of $X$, $P^0,P^1,\ldots,P^N$ with
the following properties:
\begin{enumerate}
\item For every $i\ge 1$, $P^i$ is a refinement of $P^{i-1}$.
\item $P^0=\{X\}$ and $P^N=\{\{x\}:\ x\in X\}$. \item For every
$i\ge 1$ there exists $a_i>0$ such that if $A\in P^{i-1}$ and
$B,C\in P^i$ are such that $B,C\subseteq A$, then there exists a
one-to-one onto function $\phi=\phi_{B,C}:B\to C$ such that for
every $x\in B$, $d(x,\phi(x))\le a_i$. \item
$\ell=\sqrt{\sum_{i=1}^N a_i^2}$.
\end{enumerate}
\end{definition}
For $p\ge 1$ we can can define an analogous concept if we demand
that $\ell=\left(\sum_{i=1}^N a_i^p\right)^{1/p}$. In this case we
call the parameter obtained the $\ell_p$ length of $(X,d)$, and
denote it by $\ell_p(X,d)$. Observe that it is always the case
that $\ell_p(X,d)\le \diam(X)$.

\medskip

 Recall that for $p\in [1,2]$, a Banach
space $Y$ is called $p$-smooth with constant $S$ if for every
$x,y\in Y$,
$$
\|x+y\|_Y^p+\|x-y\|_Y^p\le 2\|x\|_Y^p+2S^p\|y\|_Y^p.
$$
The least constant $S$ for which this inequality holds is called
the $p$-smoothness constant of $Y$, and is denoted $S_p(Y)$. It is
known~\cite{BCL94} that for $q\ge 2$, $S_2(L_q)\le \sqrt{q-1}$,
and for $q\in [1,2]$, $S_q(L_q)\le 1$.

The following theorem relates the notion of length to
nonembeddability results.

\begin{theorem}\label{thm:length}
Let $(X,d)$ be a metric space and $Y$ a $p$-smooth Banach space.
Then
$$
c_Y(X,d)\ge \frac{1}{2^{1-1/p}\cdot
S_p(Y)\ell_{p}(X,d)}\left(\frac{1}{|X|^2}\sum_{x,y\in
X}d(x,y)^p\right)^{1/p}.
$$
In particular for $2\le p<\infty$,
$$
c_p(X,d)\ge
\frac{1}{\ell(X,d)\sqrt{2p-2}}\left(\frac{1}{|X|^2}\sum_{x,y\in
X}d(x,y)^2\right)^{1/2}.
$$
\end{theorem}

\begin{proof} Let $\{P^i\}_{i=0}^N$, $\{a_i\}_{i=1}^N$ be as
above, and denote by ${\cal F}_i$ the $\sigma$-algebra generated
by the partition $P^i$. In what follows all expectations are taken
with respect to the uniform probability measure on $X$. Given a
bijection $f:X\to Y$ we let $f_i=\E(f|{\cal F}_i)$. In other
words, if $A\in P^i$ and $x\in A$ then
$$
f_i(x)=\frac{1}{|A|} \sum_{y\in A} f(y).
$$
Now $\{f_i\}_{i=0}^N$ is a martingale, so by Pisier's
inequality~\cite{Pisier75} (see Theorem 4.2 in~\cite{NPSS04} for
the constant we use below), we see that
$$
\E \|f_N-f_0\|^p_Y\le \frac{S_p(Y)^p}{2^{p-1}-1}\sum_{j=0}^{N-1}
\E \|f_{j+1}-f_j\|_Y^p.
$$
Now $f_0=\E f$ and $f_N=f$. Thus
\begin{eqnarray*}
\E\|f_N-f_0\|_Y^p&=&\frac{1}{|X|}\sum_{x\in
X}\left\|f(x)-\frac{1}{|X|}\sum_{y\in X} f(y)\right\|_Y^p\\&\ge&
\frac{1}{2^{p-1}|X|^2}\sum_{x,y\in X} \|f(x)-f(y)\|_Y^p\\&\ge&
\frac{1}{2^{p-1}\|f^{-1}\|_{\Lip}^p}\cdot\frac{1}{|X|^2}\sum_{x,y\in
X} d(x,y)^p.
\end{eqnarray*}
On the other hand fix $j\in \{0,\ldots, N-1\}$, and $A\in P^j$,
$B\in P^{j+1}$ such that $x\in B\subseteq A$. Then
\begin{eqnarray*}
f_{j}(x)-f_{j+1}(x)=\frac{1}{|A|} \sum_{y\in A} f(y)-\frac{1}{|B|}
\sum_{y\in B} f(y) =\frac{1}{|A|}\sum_{A\supseteq C\in
P^{j+1}}\left(\sum_{y\in C} [f(\phi_{C,B}(y))-f(y)]\right).
\end{eqnarray*}
So by convexity
$$
\|f_j(x)-f_{j+1}(x)\|_Y\le \|f||_{\Lip}\cdot a_{j+1}.
$$
It follows that
\begin{eqnarray*}
c_Y(X,d)&\ge& \frac{1}{2^{1-1/p}\cdot S_p(Y)}\cdot
\left(\frac{\frac{1}{|X|^2}\sum_{x,y\in X}d(x,y)^p}{\sum_{j=1}^N
a_j^p}\right)^{1/p}\\&=&\frac{1}{2^{1-1/p}\cdot
S_p(Y)\ell_{p}(X,d)}\left(\frac{1}{|X|^2}\sum_{x,y\in
X}d(x,y)^p\right)^{1/p}.
\end{eqnarray*}
\end{proof}

\bigskip

As shown in~\cite{MS86,Led01}, if we consider the group of
permutations of $\{1,\ldots,n\}$, $S_n$, equipped with the metric
$d(\sigma,\pi)=|\{i:\ \sigma(i)\neq \tau(i)\}|$, then
$\ell(S_n,d)\le 2\sqrt{n}$, while $\diam (S_n)=\Theta(n)$. It
follows from Theorem~\ref{thm:length} that
$c_2(S_n)=\Omega(\sqrt{n})$. On the other hand, by mapping each
permutation $\pi\in S_n$ to the matrix $({\bf 1}_{\pi(i)=j})$ we
see that $c_2(S_n)=O(\sqrt{n})$. Thus
$$
c_2(S_n)=\Theta\left(\sqrt{\frac{\log |S_n|}{\log \log
|S_n|}}\right).
$$
Similar optimal bounds can be deduced for $c_p(S_n)$, $p\ge 1$.

The metric $d$ on $S_n$ is equivalent to the shortest path metric
induced by the Cayley graph on $S_n$ obtained by taking the set of
all transpositions as generators. It is of interest to study the
Euclidean distortion of metrics on $S_n$ induced by Cayley graphs
coming from other generating sets. Recently in~\cite{Kass05} it was
shown that there exists a bounded set of generators of $S_n$ with
respect to which the Cayley graph is an expander (this settles a
long standing conjecture- see~\cite{RSW04}). It follows that there
exists a set of generators of $S_n$ with respect to which the metric
induced by the Cayley graph has Euclidean distortion $\Omega(\log
|S_n|)=\Omega(\log n\log\log n)$.

Another example discussed in~\cite{MS86,Led01} is the case of the
Hamming cube. In this case $\ell(\F_2^d,\rho)=O(\sqrt{d})$, and so
Theorem~\ref{thm:length} implies that for $p\ge 2$,
$c_p(\F_2^d,\rho)\ge c(p)\sqrt{d}$. This result was first proved
in~\cite{NS02}.

More generally, let $G$ be a finite group equipped with a
translation invariant metric $d$. Let $G=G_0\supseteq G_1\supseteq
\cdots\supseteq G_n=\{e\}$ be a decreasing sequence of subgroups.
Then it is shown in~\cite{MS86,Led01} that
$$
\ell(G,d)\le \sqrt{\sum_{j=1}^n \big[\diam(G_{i-1}/G_i)\big]^2}.
$$
This estimate implies a wide range of additional nonembeddability
results.

\section{Appendix: Quantitative estimates in Bourgain's noise sensitivity
theorem}\label{section:appendix}

In this section we prove theorem~\ref{thm:sensitive}. The proof is
a repetition of Bourgain's proof in~\cite{Bou02}, with an
optimization of the dependence on the various parameters. Since
such quantitative bounds are very useful, and they are not stated
in~\cite{Bou02}, we believe that it is worthwhile to reproduce the
argument here.

\begin{theorem}[Bourgain's distributional inequality on the
Fourier spectrum of Boolean functions]\label{thm:bourgain ineq}
Let $f:\F_2^d\to \{0,1\}$ be a Boolean function. For $2<k\le d$
and $\beta\in (0,1)$ define
$$
J_\beta=\left\{j\in \{1,\ldots,d\}:\
\sum_{\substack{A\subseteq\{1,\ldots,d\}\\|A|\le k,\ j\in A}}
\widehat f (A)^2 \ge \beta \right\}.
$$
Then
$$
 \sum_{\substack{A\subseteq\{1,\ldots,d\}\\|A|< k,\
A\setminus J_\beta\neq \emptyset}}\widehat f(A)^2\le
C^{\sqrt{\log_2(2/\delta)\log_2\log_2 k}}\cdot
\left(\delta\sqrt{k} +4^k\sqrt{\beta}\right),
$$
where
$$
\delta= \sum_{\substack{A\subseteq\{1,\ldots,d\}\\|A|\ge
k}}\widehat f(A)^2,
$$
and $C$ is a universal constant.
\end{theorem}

The quantitative version of Bourgain's noise sensitivity theorem,
as stated in Theorem~\ref{thm:sensitive}, follows from
Theorem~\ref{thm:bourgain ineq}. Indeed, we are assuming that
$f:\F_2^d\to \{-1,1\}$ satisfies
$$
\sum_{A\subseteq \{1,\ldots,d\}}(1-\e)^{|A|}\widehat f(A)^2\ge
1-\delta.
$$
By Parseval's identity, $\sum_{A\subseteq \{1,\ldots,d\}}\widehat
f(A)^2= 1$. Thus
$$
1-\delta\le \sum_{\substack{A\subseteq \{1,\ldots,d\}\\|A|<
1/\e}}\widehat f(A)^2 +\frac{1}{e}\sum_{\substack{A\subseteq
\{1,\ldots,d\}\\|A|\ge 1/\e}}\widehat f(A)^2\le 1-
\sum_{\substack{A\subseteq \{1,\ldots,d\}\\|A|\ge 1/\e}}\widehat
f(A)^2+ \frac{1}{e}\sum_{\substack{A\subseteq
\{1,\ldots,d\}\\|A|\ge 1/\e}}\widehat f(A)^2.
$$
It follows that
$$
\sum_{\substack{A\subseteq \{1,\ldots,d\}\\|A|\ge 1/\e}}\widehat
f(A)^2\le 2\delta.
$$
Now, choosing $k=1/\e$ in Theorem~\ref{thm:bourgain ineq} we get
that
$$
\beta |J_{\beta}|\le \sum_{j=1}^d
\sum_{\substack{A\subseteq\{1,\ldots,d\}\\|A|\le 1/\e,\ j\in A}}
\widehat f (A)^2=\sum_{\substack{A\subseteq \{1,\ldots,
d\}\\|A|\le 1/\e}} |A| \widehat f (A)^2\le \frac{1}{\e},
$$
i.e. $|J_\beta|\le \frac{1}{\e\beta}$. Thus, if we define
$g:\F_2^d\to \R$ by
$$
g(x)=\sum_{A\subseteq J_{1/\e}} \widehat f(A)W_A(x),
$$
then $g$ depends on at most $\frac{1}{\e\beta}$ coordinates.
Moreover, by Theorem~\ref{thm:bourgain ineq} applied to the
Boolean function $(1+f)/2$, we get that
\begin{eqnarray*}
\int_{\F_2^d}
[f(x)-g(x)]^2d\mu(x)&=&\sum_{\substack{A\subseteq\{1,\ldots,d\}\\A\setminus
J_{1/\e}\neq \emptyset}}\widehat f(A)^2\\&\le&
\sum_{\substack{A\subseteq\{1,\ldots,d\}\\|A|> 1/\e}}\widehat
f(A)^2+ \sum_{\substack{A\subseteq\{1,\ldots,d\}\\|A|\le 1/\e,\
A\setminus J_{1/\e}\neq \emptyset}}\widehat f(A)^2\\
&\le& 2\delta+4\cdot C^{\sqrt{\log_2(2/\delta)\log_2\log_2
k}}\cdot \left(\delta\sqrt{k} +4^k\sqrt{\beta}\right),
\end{eqnarray*}
as required.

\bigskip

In the proof of Theorem~\ref{thm:bourgain ineq} we will use the
following well known fact: For every $f:\F_2^d\to \R$ with
$\int_{\F_2^d}f(x)d\mu(x)=0$, and every $p\in [1,2]$,
\begin{eqnarray}\label{eq:rademacher}
\sqrt{p-1}\cdot \left(\sum_{i=1}^d \widehat
f\left(\{i\}\right)^2\right)^{1/2}\le
\left(\int_{\F_2^d}|f(x)|^pd\mu(x)\right)^{1/p}.
\end{eqnarray}
This is true since by the Bonami-Beckner
inequality~\cite{Bon70,Bec75},
\begin{eqnarray*}
\left(\int_{\F_2^d}|f(x)|^pd\mu(x)\right)^{1/p}&\ge&
\left(\int_{\F_2^d}\left|\sum_{A\subseteq\{1,\ldots,d\}}(p-1)^{\frac{|A|}{2}}\widehat
f(A)W_A(x) \right|^2\right)^{1/2}\\
&=& \left(\sum_{ A\subseteq\{1,\ldots,d\}}(p-1)^{|A|}\widehat
f(A)^2\right)^{1/2}\\
&\ge& \sqrt{p-1}\cdot \left(\sum_{i=1}^d \widehat
f\left(\{i\}\right)^2\right)^{1/2}.
\end{eqnarray*}

\begin{lemma}\label{lem:step} Fix $t,\delta,\beta\in (0,1)$ and $p\in (1,2)$.
Let $I, J$ be two disjoint finite sets and $f:\F_2^I\times
\F_2^J\to \{0,1\}$ a Boolean function. Assume that for every $i\in
I$,
$$
\sum_{\substack{A\subseteq I\cup J\\ |A|< k,\ i\in A}} \widehat
f(A)^2 \le \beta,
$$
and
$$
\sum_{\substack {A\subseteq I\cup J\\ |A|\ge k}}\widehat f(A)^2\le
\delta.
$$
Then
$$
t^{p/2} \sum_{\substack{A\subseteq I\cup J\\ |A\cap I|=1}}\widehat
f(A)^2
\le \frac{2t}{(p-1)^{p/2}}\sum_{\substack{A\subseteq I\cup J\\
|A\cap I|<k}}|A\cap I|\cdot \widehat f(A)^2 +
\frac{2\delta}{(p-1)^{p/2}}+\left(3^{k+2}\sqrt{\beta}\right)^{p/2}+(8t\delta)^{p/2}.
$$
\end{lemma}

\begin{proof} Fix $t\in (0,1)$ which will be specified later.
Let $\{s_i\}_{i\in I}$ be i.i.d. $\{0,1\}$ valued random variables
with mean $1-t$. Let $S\subseteq I$ be the random subset $S=\{i\in
I:\ s_i=1\}$. Fix $y\in \F_2^J$ and for $A\subseteq I$ denote
$$
\widehat f_y(A) =\int_{\F_2^I}
f(x,y)W_A(x)d\mu(x)=\sum_{B\subseteq J} \widehat f(A\cup B)W_B(y).
$$
Define $g_y:\F_2^I\to \R$ by
$$
g_y(x)=\sum_{\substack{A\subseteq I\\ A\not \subseteq S}}\widehat
f_y(A)W_A(x)=\sum_{\substack{A\subseteq I\\ A\not \subseteq
S}}\sum_{B\subseteq J}\widehat f(A\cup
B)W_A(x)W_B(y)=\int_{\F_2^{I\setminus S}}
f(x,y)d\mu\left((x_j)_{j\in I\setminus S}\right).
$$
Then
\begin{eqnarray}\label{eq:ident}
2\int_{\F_2^I} [g_y(x)]^2d\mu(x)=\int_{\F_2^I} |g_y(x)|d\mu(x).
\end{eqnarray}
To check this identity observe that if $\psi:\F_2^d\to \{0,1\}$ is
a Boolean function, with $\int_{\F_2^d}\psi(x)d\mu(x)=P$, then
$2\int_{\F_2^d} (\psi(x)-P)^2d\mu(x)=\int_{\F_2^d}
|\psi(x)-P|d\mu(x)=2P(1-P)$. Thus~\eqref{eq:ident} follows by
fixing $(x_j)_{j\in S}$ and $(y_j)_{j\in J}$, applying this
observation to the Boolean function $(x_j)_{j\in I\setminus
S}\mapsto f(x,y)$, and then integrating with respect to
$(x_j)_{j\in S}$.

Using~\eqref{eq:rademacher}, H\"older's inequality,
and~\eqref{eq:ident}, we get that
$$
\sqrt{p-1}\cdot \left(\sum_{i\in I\setminus S} \widehat
f_y(\{i\})^2\right)^{1/2}\le \|g_y\|_p\le
\|g_y\|_1^{\frac{2}{p}-1}\cdot
\|g_y\|_2^{\frac{2p-2}{p}}=2^{\frac{2}{p}-1}\|g_y\|_2^{\frac{2}{p}}\le
2^{\frac{2}{p}-1}\left(\sum_{\substack{A\subseteq I\\ A\not
\subseteq S}}\widehat f_y(A)^2\right)^{1/p}.
$$
In other words,
\begin{eqnarray}\label{eq:add selectors}
\left(\sum_{i\in I}(1-s_i)\widehat f_y(\{i\})^2\right)^{p/2}\le
\frac{2}{(p-1)^{p/2}}\sum_{A\subseteq I}\left(1-\prod_{i\in A}
s_i\right) \widehat f_y(A)^2.
\end{eqnarray}
Observe that $t=s_i-(1-t)+1-s_i$, so that, since $p\le 2$,
\begin{eqnarray}\label{eq:p less 2}
\nonumber\left(\sum_{i\in I}(1-s_i)\widehat
f_y(\{i\})^2\right)^{p/2}&\ge& t^{p/2} \left(\sum_{i\in I}\widehat
f_y(\{i\})^2\right)^{p/2} -\left|\sum_{i\in I}(s_i-(1-t)))\widehat
f_y(\{i\})^2\right|^{p/2}\\
&\ge& t^{p/2} \sum_{i\in I}\widehat f_y(\{i\})^2 -\left|\sum_{i\in
I}(s_i-(1-t)))\widehat f_y(\{i\})^2\right|^{p/2} .
\end{eqnarray}
Combining~\eqref{eq:add selectors} and~\eqref{eq:p less 2}, taking
expectation (with respect to $\{s_i\}_{i\in I}$), and integrating
with respect to $y\in \F_2^J$, we get that
\begin{eqnarray}\label{eq:after expect}
t^{p/2} \sum_{\substack{A\subseteq I\cup J\\ |A\cap I|=1}}\widehat
f(A)^2&\le &\frac{2}{(p-1)^{p/2}}\sum_{A\subseteq
I}\left(1-(1-t)^{|A|} \right) \widehat f(A)^2+\nonumber
\\&\phantom{\le}&\E \int_{\F_2^J} \left|\sum_{i\in
I}(s_i-(1-t)))\widehat f_y(\{i\})^2\right|^{p/2}d\mu(y) .
\end{eqnarray}

To estimate the second summand in~\eqref{eq:after expect} observe
that
\begin{eqnarray}\label{eq:steps}
\!\!\!\!\!\!\!\!\!\label{eq:jensen1}\E \int_{\F_2^J}
\left|\sum_{i\in I}(s_i-(1-t))\widehat
f_y(\{i\})^2\right|^{p/2}d\mu(y)&\le& \left(\E \int_{\F_2^J}
\left|\sum_{i\in I}(s_i-(1-t))\widehat
f_y(\{i\})^2\right|d\mu(y)\right)^{p/2}\\
&\le& \left(2\E \int_{\F_2^J} \left(\sum_{i\in I}(1-s_i)\widehat
f_y(\{i\})^4\right)^{1/2}d\mu(y)\right)^{p/2}, \label{eq:tri}
\end{eqnarray}
where in~\eqref{eq:jensen1}  we used Jensen's inequality, and
in~\eqref{eq:tri} we used the fact that if $X_1,\ldots,X_n$ are
independent mean $0$ random variables then $\E\left|X_1+\cdots
+X_n\right|\le 2\E\sqrt{X_1^2+\cdots+X_n^2}$, which can be proved
using the following standard symmetrization argument. Let
$Y_1,\ldots,Y_n$ be i.i.d. copies of $X_1,\ldots,X_n$, and let
$\e_1,\ldots,\e_n$ be i.i.d. independent $\pm 1$ Bernoulli random
variables which are also independent of the $\{X_j\}_{j=1}^n$ and
$\{Y_j\}_{j=1}^n$. Then
\begin{eqnarray*}
\E\left|\sum_{j=1}^n X_j\right|&=&\E\left|\sum_{j=1}^n (X_j-\E
Y_j)\right| \le \E\left|\sum_{j=1}^n (X_j- Y_j)\right|
= \E\left|\sum_{j=1}^n \e_j(X_j- Y_j)\right|\\
&\le& \E \sqrt{\E_\e \left|\sum_{j=1}^n \e_j(X_j- Y_j)\right|^2 }=
\E\sqrt{\sum_{j=1}^n (X_j-Y_j)^2}\le 2\E\sqrt{\sum_{j=1}^n X_j^2}.
\end{eqnarray*}

Now, from the inequality
\begin{eqnarray*}
\left|\widehat f_y(\{i\})\right|=\left|\sum_{\substack{A\subseteq I\cup J\\
A\cap I=\{i\}}} \widehat f(A)W_A(y) \right|\le \left|\sum_{\substack{A\subseteq I\cup J,\ |A|< k\\
A\cap I=\{i\}}} \widehat f(A)W_A(y) \right|+\left|\sum_{\substack{A\subseteq I\cup J,\  |A|\ge k\\
A\cap I=\{i\}}} \widehat f(A)W_A(y) \right|
\end{eqnarray*}
we get that
\begin{eqnarray*}
 \left(\sum_{i\in I}(1-s_i)\widehat
f_y(\{i\})^4\right)^{1/2}&\le& 4\left(\sum_{i\in I} \left(\sum_{\substack{A\subseteq I\cup J,\ |A|< k\\
A\cap I=\{i\}}} \widehat f(A)W_A(y)\right)^4 \right)^{1/2}+\\&\phantom{\le}& 4 \sum_{i\in I}\left(\sum_{\substack{A\subseteq I\cup J,\  |A|\ge k\\
A\cap I=\{i\}}} \widehat f(A)W_A(y) \right)^2(1-s_i).
\end{eqnarray*}
Thus, by Parseval's indentity and the Bonami-Beckner inequality we
deduce that

\begin{eqnarray}\label{eq:second beckner}
&&\!\!\!\!\!\!\!\!\!\!\!\!\!\nonumber\E \int_{\F_2^J}
\left(\sum_{i\in I}(1-s_i)\widehat
f_y(\{i\})^4\right)^{1/2}d\mu(y)\\&\le& 4\int_{\F_2^J}\left(\sum_{i\in I} \left(\sum_{\substack{A\subseteq I\cup J,\ |A|< k\\
A\cap I=\{i\}}} \widehat f(A)W_A(y)\right)^4
\right)^{1/2}d\mu(y)+\nonumber 4t\sum_{\substack {A\subseteq I\cup
J\\ |A|\ge k}}\widehat f(A)^2\\\nonumber
&\le& 4 \left(\sum_{i\in I} \int_{\F_2^J}\left(\sum_{\substack{A\subseteq I\cup J,\ |A|< k\\
A\cap I=\{i\}}} \widehat f(A)W_A(y)\right)^4d\mu(y)
\right)^{1/2}+4t\delta\\
&\le& 4\cdot 3^k\left(\sum_{i\in I} \left(\sum_{\substack{A\subseteq I\cup J,\ |A|< k\\
A\cap I=\{i\}}} \widehat f(A)^2\right)^2
\right)^{1/2}+4t\delta\nonumber\\&\le& 4\cdot 3^k\cdot \left(\max_{i\in I} \sum_{\substack{A\subseteq I\cup J,\ |A|< k\\
A\cap I=\{i\}}} \widehat f(A)^2\right)^{1/2}\cdot \left(\sum_{i\in I} \sum_{\substack{A\subseteq I\cup J,\ |A|< k\\
A\cap I=\{i\}}} \widehat f(A)^2 \right)^{1/2} +4t\delta\nonumber\\
&\le& 4\cdot 3^k \sqrt{\beta}+4t\delta.
\end{eqnarray}
Combining~\eqref{eq:second beckner} and~\eqref{eq:tri}
with~\eqref{eq:after expect}, we see that
\begin{eqnarray*}
t^{p/2} \sum_{\substack{A\subseteq I\cup J\\ |A\cap I|=1}}\widehat
f(A)^2&\le &\frac{2}{(p-1)^{p/2}}\sum_{A\subseteq
I}\left(1-(1-t)^{|A|} \right) \widehat f(A)^2+ \left(8\cdot 3^k
\sqrt{\beta}+8t\delta\right)^{p/2}\\
&\le& \frac{2t}{(p-1)^{p/2}}\sum_{\substack{A\subseteq I\cup J\\
|A\cap I|<k}}|A\cap I|\cdot \widehat f(A)^2 +
\frac{2\delta}{(p-1)^{p/2}}+\left(3^{k+2}\sqrt{\beta}\right)^{p/2}+(8t\delta)^{p/2}.
\end{eqnarray*}
\end{proof}

\begin{proof}[Proof of Theorem~\ref{thm:bourgain ineq}]
For every integer $t\ge 0$ define
$$
\rho_r=\sum_{\substack{A\subseteq \{1,\ldots,d\}\\
2^{r-1}\le |A\setminus J_\beta|<2^{r}}}\widehat f (A)^2.
$$
We also write
$$
\delta = \sum_{\substack {A\subseteq \{1,\ldots,d\}\\ |A|\ge
k}}\widehat f(A)^2.
$$

Fix $0\le r\le \log_2 k$, which will be chosen presently. Let $I$
be a uniformly random subset of $\{1,\ldots,d\}\setminus J_\beta$
of size $2^{-r-6}(d-|J_\beta|)$. If $A\subseteq \{1,\ldots, d\}$
satisfies $2^{r-1}\le |A\setminus J_\beta|<2^{r}$ then a standard
counting argument shows that $\Pr[|A\cap I|=1]\ge C$ and $\E
|A\cap I|\le C'$, where $C,C'$ are universal constants. Observe
also that by the definition of $J_\beta$, the sets $I,
J=\{1,\ldots,d\}\setminus I$ satisfy the conditions of
Lemma~\ref{lem:step}. Taking expectation with respect to $I$ of
the conclusion of lemma~\ref{lem:step}, we get that there exists a
universal constant $c>0$ such that
\begin{eqnarray}\label{eq:almost done} ct^{p/2} \rho_r \le
\frac{t}{2^r(p-1)^{p/2}}\sum_{s\le \log_2 k}
2^s\rho_s+\frac{\delta}{(p-1)^{p/2}}+\left(3^{k}\sqrt{\beta}\right)^{p/2}+(t\delta)^{p/2}.
\end{eqnarray}

Define
$$
\gamma=  \sum_{\substack{A\subseteq\{1,\ldots,d\}\\|A|< k,\
A\setminus J_\beta\neq \emptyset}}\widehat f(A)^2=\sum_{s\le
\log_2 k} \rho_s.
$$

Our goal is to prove that for some constant $C>1$,
\begin{eqnarray}\label{eq:hope}
\gamma\le C^{\sqrt{\log_2(2/\delta)\log_2\log_2
k}}\cdot\left(\delta\sqrt{k} +4^k\sqrt{\beta}\right).
\end{eqnarray}

If
$$
\sum_{s\le\log_2 k}2^s\rho_s\ge \gamma \sqrt{k}
$$
then choose $r\le \log_2 k$ to satisfy
\begin{eqnarray}\label{eq:case}
2^r\rho_r\ge \frac{1}{\log_2 k}\sum_{s\le\log_2 k}2^s\rho_s.
\end{eqnarray}
In particular it follows that $\rho_r\ge
\frac{\gamma}{\sqrt{k}\cdot \log_2k}$. Choose
$p=2-2\sqrt{\frac{\log_2\log_2k}{\log_2(1/\delta)}}$. We may
assume that $p\in (3/2,2)$, since otherwise~\eqref{eq:hope} holds
vacuously. Moreover, if
$\delta^{p/2}>\frac{c\gamma}{2\sqrt{k}\cdot \log_2k}$
then~\eqref{eq:hope} holds true. Thus, choosing
$t=\left(\frac{c}{10\log_2 k}\right)^{2/(2-p)}$
in~\eqref{eq:almost done}, and using~\eqref{eq:case}, we
obtain~\eqref{eq:hope} in this case.

It remains to deal with the case $\sum_{s\le\log_2 k}2^s\rho_s<
\gamma \sqrt{k}$. Choosing $r$ such that $\rho_r\ge
\frac{\gamma}{\log_2 k}$, and $p=1+\frac{1}{\log_2 k}$, $t\approx
\frac{1}{k(\log k)^4}$, \eqref{eq:almost done} shows
that~\eqref{eq:hope} holds true in this case as well.
\end{proof}


 \section{Acknowledgments}
 We are grateful to Keith Ball, Henry Cohn, Piotr Indyk, L\'aszl\'o Lov\'asz and Gideon Schechtman for helpful suggestions.

\bibliographystyle{abbrv}

\def\cprime{$'$}

\end{document}